\documentclass[leqno,12pt]{amsart}
\usepackage{amssymb} 
\theoremstyle{plain} 
\newtheorem{theorem}{\indent\sc Theorem}[section] 

\newtheorem{proposition}[theorem]{\indent\sc Proposition}

\theoremstyle{definition} 

\begin{document}

\title[Partial Differential system in two variables]{Partial Differential system in two variables with $W(D_6^{(1)})-$symmetry and the Garnier system in two variables \\}

\author{By\\
Yusuke Sasano}

\renewcommand{\thefootnote}{\fnsymbol{footnote}}
\footnote[0]{2000\textit{ Mathematics Subjet Classification}.
34M55; 34M45; 58F05; 32S65.}


\begin{abstract}
In this note, we will compare the Garnier system in two variables with four-dimensional partial differential system in two variables with $W(D_6^{(1)})$-symmetry. Both systems are different in each compactification in the variables $q_1,q_2$, however, has same five holomorphy conditions in the variables $p_1,p_2$.
\end{abstract}
\maketitle

\section{Motivation}
In this note, we will compare the system \eqref{eq:Gar} (which is equivalent to the Garnier system in two variables) with four-dimensional partial differential system \eqref{eq:D6} in two variables with $W(D_6^{(1)})$-symmetry. Both systems are different in each compactification in the variables $q_1,q_2$, however, has same five holomorphy conditions in the variables $p_1,p_2$.

{\small
\begin{center}
\begin{tabular}{|c||c|c|} \hline
Equation & Garnier system  \eqref{eq:Gar} & $D_6^{(1)}$ system \eqref{eq:D6}  \\ \hline
Figure & Figure 1 & Figure 3  \\ \hline
Condition of $q_1$ & $\left(\frac{1}{q_1},-(p_1 q_1+p_2 q_2+\alpha_0)q_1,\frac{q_2}{q_1},p_2 q_1 \right)$ & $\left(\frac{1}{q_1},-(p_1q_1+\alpha_2)q_1,q_2,p_2 \right)$  \\ \hline
Condition of $q_2$ & $\left(\frac{q_1}{q_2},p_1 q_2,\frac{1}{q_2},-(p_2 q_2+p_1 q_1+\alpha_0)q_2 \right)$ & $\left(q_1,p_1,\frac{1}{q_2},-(p_2q_2+\alpha_4)q_2 \right)$  \\ \hline
others & $\left(-((q_1-t)p_1-\alpha_5)p_1,\frac{1}{p_1},q_2,p_2 \right)$ & $\left(-((q_1-t)p_1-\alpha_0)p_1,\frac{1}{p_1},q_2,p_2 \right)$  \\ \hline
 & $\left(-((q_1-s)p_1-\alpha_4)p_1,\frac{1}{p_1},q_2,p_2 \right)$ & $\left(-((q_1-s)p_1-\alpha_1)p_1,\frac{1}{p_1},q_2,p_2 \right)$  \\ \hline
 & $\left(-\left( (q_1-q_2)p_1-\alpha_3 \right)p_1,\frac{1}{p_1},q_2,p_2+p_1 \right)$ & $\left(-\left( (q_1-q_2)p_1-\alpha_3 \right)p_1,\frac{1}{p_1},q_2,p_2+p_1 \right)$  \\ \hline
 & $\left(q_1,p_1,-((q_2-1)p_2-\alpha_1)p_2,\frac{1}{p_2} \right)$ & $\left(q_1,p_1,-((q_2-1)p_2-\alpha_5)p_2,\frac{1}{p_2} \right)$  \\ \hline
 & $\left(q_1,p_1,-(q_2 p_2-\alpha_2)p_2,\frac{1}{p_2} \right)$ & $\left(q_1,p_1,-(q_2p_2-\alpha_6)p_2,\frac{1}{p_2} \right)$  \\ \hline
Holomorphy & \eqref{holo;SGar} & \eqref{holo;D6}  \\ \hline
Symmetry & \eqref{symmetry:SGar} & \eqref{symmetry:D6}  \\ \hline
\end{tabular}
\end{center}
}

We remark that each compactification can be obtained by successive blowing-ups and blowing-downs (containing ${\mathbb P}^2$-flop) in 4-dimensional projective space ${\mathbb P}^4$, respectively (see \cite{sasa3,sasa7}).

At first, in August in 2015 we constructed the system \eqref{eq:D6} by modifying the holomorphy conditions (cf. \cite{sasa3}) of 2-coupled Painlev\'e VI system of type $D_6^{(1)}$. In those days, this system is doubted as padding out another time-variable $s$ (cf. \cite{Fuji}). For the system \eqref{eq:D6} we showed its Completely integrable, Holomorphy conditions and Affine Weyl group symmetries.

Next, we constructed the system  \eqref{eq:Gar} by modifying the holomorphy conditions \eqref{holo;D6} (cf. \cite{sasa3}) of $D_6^{(1)}$ system \eqref{eq:D6}. By using the holomorphy condition
\begin{align}
\begin{split}
\left(\frac{q_1}{q_2},p_1 q_2,\frac{1}{q_2},-(p_2 q_2+p_1 q_1+\alpha_0)q_2 \right),
\end{split}
\end{align}
we can make a change of variables \eqref{symp} for the system \eqref{eq:Gar} since the boundary conditions in the variables $p_1,p_2$ for the system \eqref{eq:Gar} is different from one of the Garnier system \eqref{canGar} in two variables. Thanks to this transformation, we see that the system \eqref{eq:Gar} is equivalent to the Garnier system \eqref{canGar} in two variables. 

It is still an open question whether from a viewpoint of higher-dimensional minimal model  both systems  \eqref{canGar}, \eqref{eq:D6} are different or not.

Moreover, it is also still an open question whether from the viewpoint of Lax pair both systems  \eqref{canGar}, \eqref{eq:D6} are different or not (cf. \cite{Fuji2}).

\section{Garnier system in two variables}

{\bf Holomorphy conditions}

Define birational and symplectic transformations $r_i \ (i=0,1,\ldots,5)$ as follows:

\begin{align}\label{holo;SGar}
\begin{split}
r_0:(x_0,y_0,z_0,w_0)=&\left(\frac{1}{q_1},-(p_1 q_1+p_2 q_2+\alpha_0)q_1,\frac{q_2}{q_1},p_2 q_1 \right),\\
r_1:(x_1,y_1,z_1,w_1)=&\left(q_1,p_1,-((q_2-1)p_2-\alpha_1)p_2,\frac{1}{p_2} \right),\\
r_2:(x_2,y_2,z_2,w_2)=&\left(q_1,p_1,-(q_2 p_2-\alpha_2)p_2,\frac{1}{p_2} \right),\\
r_3:(x_3,y_3,z_3,w_3)=&\left(-\left( (q_1-q_2)p_1-\alpha_3 \right)p_1,\frac{1}{p_1},q_2,p_2+p_1 \right),\\
r_4:(x_4,y_4,z_4,w_4)=&\left(-((q_1-s)p_1-\alpha_4)p_1,\frac{1}{p_1},q_2,p_2 \right),\\
r_5:(x_5,y_5,z_5,w_5)=&\left(-((q_1-t)p_1-\alpha_5)p_1,\frac{1}{p_1},q_2,p_2 \right).
\end{split}
\end{align}
There exist two polynomials $H_1$ and $H_2$, such that the Hamiltonian system
\begin{equation}\label{eq:Gar}
  \left\{
  \begin{aligned}
   dq_1 =&\frac{\partial H_1}{\partial p_1}dt+\frac{\partial H_2}{\partial p_1}ds, \quad dp_1 =-\frac{\partial H_1}{\partial q_1}dt-\frac{\partial H_2}{\partial q_1}ds,\\
   dq_2 =&\frac{\partial H_1}{\partial p_2}dt+\frac{\partial H_2}{\partial p_2}ds, \quad dp_2 =-\frac{\partial H_1}{\partial q_2}dt-\frac{\partial H_2}{\partial q_2}ds
   \end{aligned}
  \right. 
\end{equation}
is transformed into a polynomial Hamiltonian system under the action of each $r_i \ (i=0,1,\ldots,5)$, where two polynomial Hamiltonians $H_1,H_2$ are given by
\begin{align}
\begin{split}
&(t - s) (t - 1) t (2 \alpha_0 + \alpha_1 + \alpha_2 + \alpha_3 + \alpha_4 + \alpha_5) 
  H_1\\
&= (t - 1) \{ (q_1 - s) (q_1 - t) q_1 p_1^2 + \alpha_0 (\alpha_0 + \alpha_1) 
       q_1 + p_1 ((t - 1) (q_1 (\alpha_3 + \alpha_4) - 
            s \alpha_3)\\
& + (2 \alpha_0 + \alpha_1) (q_1 - 
            s) q_1 + \alpha_2 (q_1 - s) (q_1 - 1) + \alpha_3 (q_1 - 
            s) - \alpha_4 (s - 1) q_1 ) \} \\
& - s \{ (q_2 - t) (q_2 - 1) q_2 p_2^2 + \alpha_0 (\alpha_0 + \alpha_4) q_2 + 
      p_2 (2 \alpha_0 (q_2 - 1) q_2 + \alpha_1 (t - 
            1) q_2\\
& + \alpha_2 t (q_2 - 1) + \alpha_4 (q_2 - 1) q_2) \} \\
& + (q_1 - s) (q_1 - t)^2 p_1^2 q_2 + 
     q_1 p_2^2 (q_2 - t) (q_2 - 1) q_2 + \alpha_0 ( \alpha_0 q_1 q_2 + \alpha_2 t q_1 + \alpha_4 t q_2 )\\
& + p_1 (2 p_2 (q_1 - s) (q_1 - t) (q_2 - 1) q_2 + 
        2 \alpha_0 (q_1 - s) (q_1 - t) q_2 + \alpha_2 (q_1 - s) (q_1 - 1) \\
& + \alpha_4 (t - s) (q_1 - t) q_2 ) + 
     p_2 (2 \alpha_0 q_1 (q_2 - 1) q_2 + \alpha_1 (t - 1) q_1 q_2 + \alpha_2 t q_1 (q_2 - 1) + \alpha_4 t (q_2 - 1) q_2), \\
H_2 &=\pi(H_1), \quad \pi=\{ t \leftrightarrow s, \ \alpha_4 \leftrightarrow \alpha_5 \}.
\end{split}
\end{align}
We note that the holomorphy conditions should be read that in the Hamiltonian $H_1$
\begin{align*}
\begin{split}
&r_5(H_1 - p_1)
\end{split}
\end{align*}
are polynomials with respect to $x_5,y_5,z_5,w_5$, and in the Hamiltonian $H_2$
\begin{align*}
\begin{split}
&r_4(H_2 -p_1)
\end{split}
\end{align*}
are polynomials with respect to $x_4,y_4,z_4,w_4$.

\vspace{0.5cm}

{\bf Completely integrable}

(see Definition; page 134: On the polynomial Hamiltonian structures of the Garnier systems; H. Kimura and K. Okamoto)

\begin{proposition}
Setting
\begin{equation}
K_1:=-H_1+\frac{ \alpha_0 \alpha_2 ( \rm{Log \rm}(s-t)-\rm{Log \rm}(s-1) )}{(2\alpha_0+\alpha_1+\alpha_2+\alpha_3+\alpha_4+\alpha_5)(t-1)^2}, \quad K_2:=-H_2.
\end{equation}
Two Hamiltonians $K_1$ and $K_2$ satisfy
\begin{equation}
\{K_1,K_2\}+\left(\frac{\partial}{\partial s} \right)K_1-\left(\frac{\partial}{\partial t} \right)K_2=0,
\end{equation}
where  $\{,\}$ denotes the Poisson brackets:
\begin{equation}
\{L_1,L_2\}=\frac{\partial L_1}{\partial p_1}\frac{\partial L_2}{\partial q_1}-\frac{\partial L_1}{\partial q_1}\frac{\partial L_2}{\partial p_1}+\frac{\partial L_1}{\partial p_2}\frac{\partial L_2}{\partial q_2}-\frac{\partial L_1}{\partial q_2}\frac{\partial L_2}{\partial p_2}.
\end{equation}
\end{proposition}
We note that the system \eqref{eq:Gar} is invariant under the birational and symplectic transformations $s_0,s_1,\ldots,s_8$ defined as follows$:$ with {\it the notation} $(*):=(q_1,p_1,q_2,p_2,t,s;\alpha_0,\alpha_1,\ldots,\alpha_5)$\rm{: \rm}
{\small
\begin{align}\label{symmetry:SGar}
\begin{split}
s_0:(*) \rightarrow &\left(\frac{1}{q_1},-(p_1 q_1+p_2 q_2+\alpha_0)q_1,\frac{q_2}{q_1},p_2 q_1,\frac{1}{t},\frac{1}{s};\alpha_0,\alpha_3,\alpha_2,\alpha_1,\alpha_4,\alpha_5 \right),\\
s_1:(*) \rightarrow &\left(q_1,p_1,q_2,p_2-\frac{\alpha_1}{q_2-1},t,s;\alpha_0+\alpha_1,-\alpha_1,\alpha_2,\alpha_3,\alpha_4,\alpha_5 \right),\\
s_2:(*) \rightarrow &\left(q_1,p_1,q_2,p_2-\frac{\alpha_2}{q_2},t,s;\alpha_0+\alpha_2,\alpha_1,-\alpha_2,\alpha_3,\alpha_4,\alpha_5 \right),\\
s_3:(*) \rightarrow &\left(q_1,p_1-\frac{\alpha_3}{q_1-q_2},q_2,p_2+\frac{\alpha_3}{q_1-q_2},t,s;\alpha_0+\alpha_3,\alpha_1,\alpha_2,-\alpha_3,\alpha_4,\alpha_5 \right),\\
s_4:(*) \rightarrow &\left(q_1,p_1-\frac{\alpha_4}{q_1-s},t,s;;\alpha_0+\alpha_4,\alpha_1,\alpha_2,\alpha_3,-\alpha_4,\alpha_5 \right),\\
s_5:(*) \rightarrow &\left(q_1,p_1-\frac{\alpha_5}{q_1-t},t,s;\alpha_0+\alpha_5,\alpha_1,\alpha_2,\alpha_3,\alpha_4,-\alpha_5 \right),\\
        s_6: (*) \rightarrow  &(q_1,p_1,q_2,p_2,s,t;\alpha_0,\alpha_1,\alpha_2,\alpha_3,\alpha_5,\alpha_4 ), \\
        s_7: (*) \rightarrow  &(1-q_1,-p_1,1-q_2,-p_2,1-t,1-s;\alpha_0,\alpha_2,\alpha_1,\alpha_3,\alpha_4,\alpha_5), \\
        s_8: (*) \rightarrow  &(\frac{(s-t)(-q_1+s q_2)}{(s-1)(t-q_1+s q_2-t q_2)},-\frac{(s-1)(t-q_1+s q_2-t q_2)(t p_1-q_1 p_1-q_2 p_2-\alpha_0)}{(s-t)t},\\
&\frac{(t-s)q_2}{-t+q_1-s q_2+t q_2},\frac{(t-q_1+s q_2-t q_2)(p_2+s p_1-q_1 p_1-q_2 p_2-\alpha_0)}{s-t},\frac{s-t}{s-1},\frac{s}{s-1};\\
&\alpha_0,\alpha_5,\alpha_2,\alpha_3,\alpha_4,\alpha_1 ).
\end{split}
\end{align}
}
The group $<s_0,s_6,s_7,s_8>$ is isomorphic to symmetric group of degree five.

The transformation $s_0$ can be exchanged the constant parameters $\alpha_1$ and $\alpha_3$. The diagonal accessible singular curve $C_3( \cong {\mathbb P}^1)$ and the accessible singular curve $C_1( \cong {\mathbb P}^1)$ (see \eqref{accessible sing.}, Figure 1.) can be transformed by the transformation $s_0$.

We remark that thanks to the transformation \eqref{symp}, we can obtain well-known Tsuda's transformation (see \cite{Tsuda1,Tsuda2,Tsuda3,Tsuda4}).

We note that the system \eqref{eq:Gar} has the following invariant divisors\rm{:\rm}
\begin{center}
\begin{tabular}{|c|c|c|} \hline
parameter's relation & invariant divisor \\ \hline
$\alpha_0=0$ & $f_0^{(1)}:=p_1, \ f_0^{(2)}:=p_2$  \\ \hline
$\alpha_1=0$ & $f_1:=q_2-1$  \\ \hline
$\alpha_2=0$ & $f_2:=q_2$  \\ \hline
$\alpha_3=0$ & $f_3:=q_1-q_2$  \\ \hline
$\alpha_4=0$ & $f_4:=q_1-s$  \\ \hline
$\alpha_5=0$ & $f_5:=q_1-t$  \\ \hline
\end{tabular}
\end{center}

\begin{figure}
\unitlength 0.1in
\begin{picture}( 41.1000, 21.6000)( 13.1000,-23.5000)
%
\special{pn 8}%
\special{pa 2430 888}%
\special{pa 1758 2190}%
\special{fp}%
\special{pa 1758 2190}%
\special{pa 5164 2190}%
\special{fp}%
%
\special{pn 8}%
\special{pa 2424 888}%
\special{pa 4432 888}%
\special{fp}%
%
\special{pn 8}%
\special{pa 4418 888}%
\special{pa 5158 2190}%
\special{fp}%
%
\special{pn 8}%
\special{pa 1758 2182}%
\special{pa 2604 1834}%
\special{fp}%
%
\special{pn 8}%
\special{pa 2610 1826}%
\special{pa 4324 1826}%
\special{dt 0.045}%
\put(30.1100,-21.0100){\makebox(0,0)[lb]{$C_3$}}%
\put(21.1600,-17.6200){\makebox(0,0)[lb]{$C_0$}}%
\put(45.0500,-15.7500){\makebox(0,0)[lb]{$C_5$}}%
%
\special{pn 20}%
\special{pa 4418 888}%
\special{pa 5158 2198}%
\special{fp}%
\put(47.0400,-13.9700){\makebox(0,0)[lb]{$C_6$}}%
\put(22.5000,-14.3000){\makebox(0,0)[lb]{$C_1$}}%
%
\special{pn 8}%
\special{pa 2438 896}%
\special{pa 3504 410}%
\special{dt 0.045}%
%
\special{pn 8}%
\special{pa 3518 434}%
\special{pa 4444 888}%
\special{dt 0.045}%
%
\special{pn 8}%
\special{pa 3490 418}%
\special{pa 3490 1502}%
\special{dt 0.045}%
%
\special{pn 8}%
\special{pa 2618 1818}%
\special{pa 3498 1510}%
\special{dt 0.045}%
%
\special{pn 8}%
\special{pa 3490 1510}%
\special{pa 4344 1826}%
\special{dt 0.045}%
%
\special{pn 20}%
\special{pa 3172 2182}%
\special{pa 3458 1326}%
\special{fp}%
\put(16.8000,-19.6000){\makebox(0,0)[lb]{$r_5$}}%
\put(18.2000,-14.7000){\makebox(0,0)[lb]{$r_4$}}%
\put(29.9000,-23.8000){\makebox(0,0)[lb]{$r_3$}}%
\put(46.6000,-18.6000){\makebox(0,0)[lb]{$r_1$}}%
\put(50.7000,-20.6000){\makebox(0,0)[lb]{$r_2$}}%
\put(32.7000,-11.7000){\makebox(0,0)[lb]{$r_0$}}%
%
\special{pn 8}%
\special{pa 2440 890}%
\special{pa 2620 1830}%
\special{dt 0.045}%
%
\special{pn 8}%
\special{pa 4420 900}%
\special{pa 4340 1830}%
\special{dt 0.045}%
%
\special{pn 8}%
\special{pa 4350 1830}%
\special{pa 5160 2200}%
\special{dt 0.045}%
%
\special{pn 20}%
\special{pa 2140 1480}%
\special{pa 2540 1480}%
\special{fp}%
%
\special{pn 20}%
\special{pa 4370 1430}%
\special{pa 4820 2050}%
\special{fp}%
%
\special{pn 20}%
\special{pa 1900 1930}%
\special{pa 2590 1670}%
\special{fp}%
%
\special{pn 8}%
\special{pa 3500 420}%
\special{pa 3690 320}%
\special{fp}%
\special{sh 1}%
\special{pa 3690 320}%
\special{pa 3622 334}%
\special{pa 3644 346}%
\special{pa 3640 370}%
\special{pa 3690 320}%
\special{fp}%
%
\special{pn 8}%
\special{pa 3500 1510}%
\special{pa 3760 1380}%
\special{fp}%
\special{sh 1}%
\special{pa 3760 1380}%
\special{pa 3692 1392}%
\special{pa 3712 1404}%
\special{pa 3710 1428}%
\special{pa 3760 1380}%
\special{fp}%
\put(34.7000,-3.6000){\makebox(0,0)[lb]{$\frac{1}{q_1}$}}%
\put(35.0000,-14.0000){\makebox(0,0)[lb]{$\frac{1}{q_2}$}}%
%
\special{pn 8}%
\special{pa 1760 2180}%
\special{pa 1570 2340}%
\special{fp}%
\special{sh 1}%
\special{pa 1570 2340}%
\special{pa 1634 2312}%
\special{pa 1612 2306}%
\special{pa 1608 2282}%
\special{pa 1570 2340}%
\special{fp}%
%
\special{pn 8}%
\special{pa 5160 2190}%
\special{pa 5400 2350}%
\special{fp}%
\special{sh 1}%
\special{pa 5400 2350}%
\special{pa 5356 2296}%
\special{pa 5356 2320}%
\special{pa 5334 2330}%
\special{pa 5400 2350}%
\special{fp}%
\put(13.1000,-23.4000){\makebox(0,0)[lb]{$\frac{1}{p_1}$}}%
\put(54.2000,-23.5000){\makebox(0,0)[lb]{$\frac{1}{p_2}$}}%
\end{picture}%
\label{CPD6fig5}
\caption{This figure denotes the boundary divisor ${\mathcal H}^{Kim}$ of ${\mathcal S}^{Kim}$, where the symbol ${\mathcal S}^{Kim}$ denotes a Hirzebruch manifold constructed by K. Kimura (see \cite{K,KOka}). The bold lines $C_i \ (i=0,1,3,5,6)$ (see \eqref{accessible sing.}) in ${\mathcal H}^{Kim}$ denote the accessible singular loci of the system \eqref{eq:Gar}. }
\end{figure}
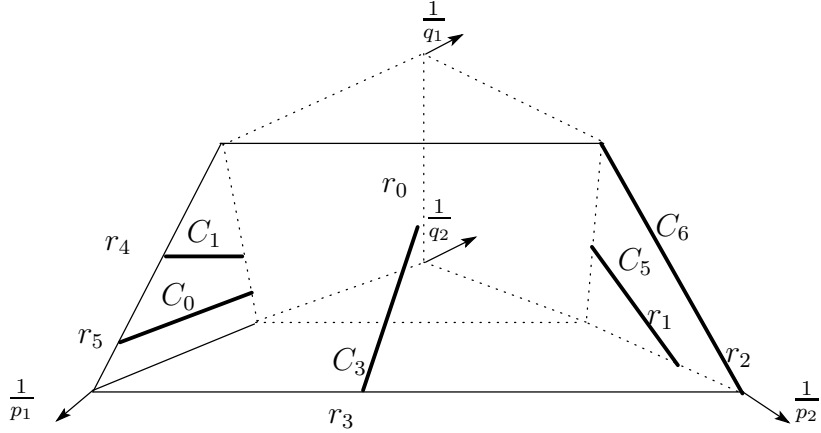

\begin{proposition}
The birational and symplectic transformation $\varphi$\rm{:\rm}
\begin{equation}\label{symp}
  \left\{
  \begin{aligned}
   Q_1=&\frac{q_1-q_2}{(s-1)q_2},\\
   P_1=&(s-1)p_1 q_2,\\
   Q_2=&\frac{s(q_2-1)}{(s-1)q_2},\\
   P_2=&\frac{(s-1)q_2(q_1p_1+q_2p_2+\alpha_0)}{s},\\
   T=&\frac{t-1}{s-1},\\
   S=&\frac{s(t-1)}{t(s-1)}
   \end{aligned}
  \right. 
\end{equation}
takes the system \eqref{eq:Gar} into the Garnier system in two variables with canonical polynomial Hamiltonians $L_1,L_2;$ {\rm (see \cite{Gar1,K,KOka,oka4,KNS,Yamada}) \rm}
\begin{align}\label{canGar}
\begin{split}
dq_1&=\frac{\partial L_1}{\partial p_1}dt+\frac{\partial L_2}{\partial p_1}ds, \quad dp_1=-\frac{\partial L_1}{\partial q_1}dt-\frac{\partial L_2}{\partial q_1}ds,\\
dq_2&=\frac{\partial L_1}{\partial p_2}dt+\frac{\partial L_2}{\partial p_2}ds, \quad dp_2=-\frac{\partial L_1}{\partial q_2}dt-\frac{\partial L_2}{\partial q_2}ds,\\
L_1 &=H_{VI}(q_1,p_1,t;\alpha_1+\alpha_5,\alpha_2,\alpha_0,\alpha_4,\alpha_3)\\
&-\frac{\alpha_1 s}{t(t-s)}q_1p_1+\frac{(s-1)q_2p_2(2q_1p_1-\alpha_3)}{(t-s)(t-1)}\\
&-\frac{t(q_1p_1-\alpha_3)p_1q_2+s(q_2p_2-\alpha_1)q_1p_2}{t(t-s)}+\frac{\{2(q_1p_1+\alpha_0)+q_2p_2+\alpha_2\}q_1q_2p_2}{t(t-1)},\\
L_2&=\pi(L_1),
\end{split}
\end{align}
where the transformation $\pi$ is explicitly given by
\begin{align}
\begin{split}
\pi:&(q_1,p_1,q_2,p_2,t,s;\alpha_0,\alpha_1,\alpha_2,\alpha_3,\alpha_4,\alpha_5) \rightarrow(q_2,p_2,q_1,p_1,s,t;\alpha_0,\alpha_3,\alpha_2,\alpha_1,\alpha_4,\alpha_5),
\end{split}
\end{align}
where for notational convenience, we have renamed $Q_i,P_i,T,S$ to $q_i,p_i,t,s$ (which are not the same as the previous $q_i,p_i,t,s$).
\end{proposition}

\begin{figure}
\unitlength 0.1in
\begin{picture}( 52.3200, 33.1200)( 10.2900,-34.6800)
%
\special{pn 8}%
\special{pa 3828 356}%
\special{pa 1866 1402}%
\special{dt 0.045}%
\special{pa 1886 1388}%
\special{pa 5444 1388}%
\special{dt 0.045}%
\special{pa 3838 356}%
\special{pa 5424 1388}%
\special{dt 0.045}%
%
\special{pn 8}%
\special{pa 1906 1388}%
\special{pa 1276 3226}%
\special{dt 0.045}%
\special{pa 1276 3226}%
\special{pa 6174 3226}%
\special{dt 0.045}%
\special{pa 5424 1394}%
\special{pa 6174 3226}%
\special{dt 0.045}%
%
\special{pn 8}%
\special{pa 3838 380}%
\special{pa 3828 2024}%
\special{dt 0.045}%
\special{pa 1906 1378}%
\special{pa 2646 2470}%
\special{dt 0.045}%
\special{pa 2646 2470}%
\special{pa 4804 2470}%
\special{dt 0.045}%
\special{pa 5424 1394}%
\special{pa 4824 2478}%
\special{dt 0.045}%
\special{pa 2646 2478}%
\special{pa 3828 2024}%
\special{dt 0.045}%
\special{pa 3828 2024}%
\special{pa 4814 2462}%
\special{dt 0.045}%
\special{pa 1286 3226}%
\special{pa 2656 2478}%
\special{dt 0.045}%
\special{pa 4832 2470}%
\special{pa 6174 3226}%
\special{dt 0.045}%
%
\special{pn 20}%
\special{pa 3848 356}%
\special{pa 4104 266}%
\special{fp}%
\special{sh 1}%
\special{pa 4104 266}%
\special{pa 4034 268}%
\special{pa 4054 284}%
\special{pa 4048 306}%
\special{pa 4104 266}%
\special{fp}%
%
\special{pn 20}%
\special{pa 1286 3226}%
\special{pa 1198 3468}%
\special{fp}%
\special{sh 1}%
\special{pa 1198 3468}%
\special{pa 1240 3412}%
\special{pa 1216 3418}%
\special{pa 1202 3400}%
\special{pa 1198 3468}%
\special{fp}%
%
\special{pn 20}%
\special{pa 6174 3226}%
\special{pa 6242 3468}%
\special{fp}%
\special{sh 1}%
\special{pa 6242 3468}%
\special{pa 6242 3398}%
\special{pa 6228 3418}%
\special{pa 6204 3410}%
\special{pa 6242 3468}%
\special{fp}%
%
\special{pn 20}%
\special{pa 3838 2016}%
\special{pa 4124 1910}%
\special{fp}%
\special{sh 1}%
\special{pa 4124 1910}%
\special{pa 4054 1914}%
\special{pa 4074 1928}%
\special{pa 4068 1952}%
\special{pa 4124 1910}%
\special{fp}%
\put(41.7300,-3.2600){\makebox(0,0)[lb]{$\frac{1}{q_1}$}}%
\put(40.3000,-17.8000){\makebox(0,0)[lb]{$\frac{1}{q_2}$}}%
\put(62.6100,-35.3700){\makebox(0,0)[lb]{$\frac{1}{p_2}$}}%
\put(10.2900,-35.9000){\makebox(0,0)[lb]{$\frac{1}{p_1}$}}%
\put(35.9100,-11.1400){\makebox(0,0)[lb]{$r_0$}}%
%
\special{pn 20}%
\special{pa 1300 3210}%
\special{pa 2660 2490}%
\special{fp}%
%
\special{pn 20}%
\special{pa 1920 1390}%
\special{pa 2670 2470}%
\special{fp}%
%
\special{pn 20}%
\special{pa 3160 3220}%
\special{pa 3020 1860}%
\special{fp}%
%
\special{pn 20}%
\special{pa 3540 1790}%
\special{pa 3990 3240}%
\special{fp}%
%
\special{pn 20}%
\special{pa 5420 1400}%
\special{pa 6170 3240}%
\special{fp}%
\put(15.6000,-28.9000){\makebox(0,0)[lb]{$r_3$}}%
\put(20.6000,-18.7000){\makebox(0,0)[lb]{$r_2$}}%
\put(29.2000,-28.0000){\makebox(0,0)[lb]{$r_5$}}%
\put(36.3000,-28.1000){\makebox(0,0)[lb]{$r_4$}}%
\put(56.1000,-24.5000){\makebox(0,0)[lb]{$r_1$}}%
\end{picture}%
\label{'SŽŸŒ³‹óŠÔGarniKimura2}
\caption{Hirzebruch manifold defined by H. Kimura and some accessible singularities}
\end{figure}
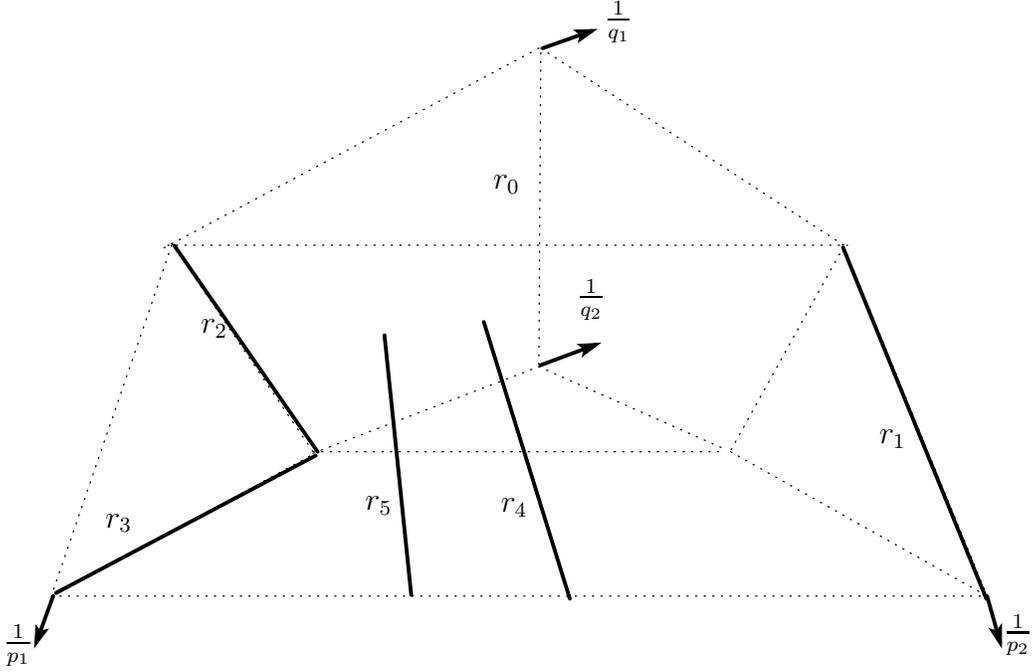

Here, the symbol $H_{VI}$ denotes
\begin{align}\label{SPVI}
\begin{split}
&H_{VI}(q,p,t;\beta_0,\beta_1,\beta_2,\beta_3,\beta_4)=\frac{1}{t(t-1)}[p^2(q-t)(q-1)q-\{(\beta_0-1)(q-1)q+\beta_3(q-t)q\\
&+\beta_4(q-t)(q-1)\}p+\beta_2(\beta_1+\beta_2)(q-t)] \quad (\beta_0+\beta_1+2\beta_2+\beta_3+\beta_4=1).
\end{split}
\end{align}
We note that this system admits the following holomorphy conditions;
{\small
\begin{align}
\begin{split}
&r_0:x_0=\frac{1}{q_1},\ y_0=-(q_1p_1+q_2p_2+\alpha_0)q_1, \ z_0=\frac{q_2}{q_1}, \  w_0=q_1p_2,\\
&r_1:x_1=q_1, \ y_1=p_1, \ z_1=-(q_2 p_2-\alpha_1)p_2, \ w_1=\frac{1}{p_2},\\
&r_2:x_2=\frac{1}{q_1},\ y_2=-(q_1p_1+q_2p_2+\alpha_2+\alpha_0)q_1, \ z_2=\frac{q_2}{q_1}, \  w_2=q_1p_2,\\
&r_3:x_3=-(q_1p_1-\alpha_3)p_1,\ y_3=\frac{1}{p_1},\ z_3=q_2,\ w_3=p_2, \\
&r_4:x_4=-((q_1+q_2-1)p_1-\alpha_4)p_1, \ y_4=\frac{1}{p_1}, \ z_4=q_2, \ w_4=p_2-p_1, \\
&r_5:x_5=-((q_1+tq_2/s-t)p_1-\alpha_5)p_1, \ y_5=\frac{1}{p_1}, \  z_5=q_2, \ w_5=p_2-\frac{tp_1}{s}.
\end{split}
\end{align}
}

\section{$D_6^{(1)}$ system}

{\bf Holomorphy conditions}

Define birational and symplectic transformations $R_i \ (i=0,1,\ldots,6)$ as follows:

\begin{align}\label{holo;D6}
\begin{split}
R_0:(x_0,y_0,z_0,w_0)=&\left(-((q_1-t)p_1-\alpha_0)p_1,\frac{1}{p_1},q_2,p_2 \right),\\
R_1:(x_1,y_1,z_1,w_1)=&\left(-((q_1-s)p_1-\alpha_1)p_1,\frac{1}{p_1},q_2,p_2 \right),\\
R_2:(x_2,y_2,z_2,w_2)=&\left(\frac{1}{q_1},-(p_1q_1+\alpha_2)q_1,q_2,p_2 \right),\\
R_3:(x_3,y_3,z_3,w_3)=&\left(-\left( (q_1-q_2)p_1-\alpha_3 \right)p_1,\frac{1}{p_1},q_2,p_2+p_1 \right),\\
R_4:(x_4,y_4,z_4,w_4)=&\left(q_1,p_1,\frac{1}{q_2},-(p_2q_2+\alpha_4)q_2 \right),\\
R_5:(x_5,y_5,z_5,w_5)=&\left(q_1,p_1,-((q_2-1)p_2-\alpha_5)p_2,\frac{1}{p_2} \right),\\
R_6:(x_6,y_6,z_6,w_6)=&\left(q_1,p_1,-(q_2p_2-\alpha_6)p_2,\frac{1}{p_2} \right).
\end{split}
\end{align}
There exist two polynomials $H_1^{D_6^{(1)}}$ and $H_2^{D_6^{(1)}}$, such that the Hamiltonian system
\begin{equation}\label{eq:D6}
  \left\{
  \begin{aligned}
   dq_1 =&\frac{\partial H_1^{D_6^{(1)}}}{\partial p_1}dt+\frac{\partial H_2^{D_6^{(1)}}}{\partial p_1}ds, \quad dp_1 =-\frac{\partial H_1^{D_6^{(1)}}}{\partial q_1}dt-\frac{\partial H_2^{D_6^{(1)}}}{\partial q_1}ds,\\
   dq_2 =&\frac{\partial H_1^{D_6^{(1)}}}{\partial p_2}dt+\frac{\partial H_2^{D_6^{(1)}}}{\partial p_2}ds, \quad dp_2 =-\frac{\partial H_1^{D_6^{(1)}}}{\partial q_2}dt-\frac{\partial H_2^{D_6^{(1)}}}{\partial q_2}ds
   \end{aligned}
  \right. 
\end{equation}
is transformed into a polynomial Hamiltonian system under the action of each $R_i \ (i=0,1,\ldots,6)$, where two polynomial Hamiltonians $H_1^{D_6^{(1)}},H_2^{D_6^{(1)}}$ are given by {\rm (cf. \cite{KNS,Fuji,sasa3}) \rm}
\begin{align}
\begin{split}
H_1^{D_6^{(1)}} &=H_{VI}(q_1,p_1,t,s;\alpha_0,\alpha_1,\alpha_2,\alpha_3+2\alpha_4+\alpha_5,\alpha_3+\alpha_6)\\
&+H_{VI}(q_2,p_2,t,s;\alpha_0+2\alpha_2+\alpha_3,\alpha_1+\alpha_3,\alpha_4,\alpha_5,\alpha_6)\\
&+\frac{2(q_1-s)q_2 \{(q_1-t)p_1+\alpha_2 \} \{(q_2-1)p_2+\alpha_4 \}}{(\alpha_0+\alpha_1+2\alpha_2+2\alpha_3+2\alpha_4+\alpha_5+\alpha_6)t(t-1)(t-s)},\\
H_2^{D_6^{(1)}} &=\pi(H_1^{D_6^{(1)}}), \quad \pi=\{ t \leftrightarrow s, \ \alpha_0 \leftrightarrow \alpha_1 \}.
\end{split}
\end{align}
The symbol $H_{VI}(q,p,t,\eta;\beta_0,\beta_1,\beta_2,\beta_3,\beta_4)$ denotes (see \cite{sasa3})
\begin{align}
\begin{split}
&t(t-1)(t-\eta) H_{VI}(q,p,t,\eta;\beta_0,\beta_1,\beta_2,\beta_3,\beta_4)\\
&=q(q-1)(q-\eta)(q-t)p^2+\{\beta_1(t-\eta)q(q-1)+2\beta_2 q(q-1)(q-\eta)\\
&+\beta_3 (t-1)q(q-\eta)+\beta_4 t(q-1)(q-\eta) \}p\\
&+\beta_2 \{(\beta_1+\beta_2)(t-\eta)+\beta_2(q-1)+\beta_3(t-1)+t \beta_4 \}q, \quad (\beta_0+\beta_1+2\beta_2+\beta_3+\beta_4=1).
\end{split}
\end{align}

\vspace{0.5cm}

We note that the holomorphy conditions should be read that in the Hamiltonian $H_1$
\begin{align*}
\begin{split}
&R_0(H_1 - p_1)
\end{split}
\end{align*}
are polynomials with respect to $x_0,y_0,z_0,w_0$, and in the Hamiltonian $H_2$
\begin{align*}
\begin{split}
&R_1(H_2 -p_1)
\end{split}
\end{align*}
are polynomials with respect to $x_1,y_1,z_1,w_1$.

\vspace{0.5cm}

{\bf Completely integrable}

\begin{proposition}
Setting
\begin{equation}
K_1:=-H_1+\frac{(\alpha_2 \alpha_3+\alpha_2 \alpha_6+\alpha_4 \alpha_6)( \rm{Log \rm}(s-t)-\rm{Log \rm}(s-1) )}{(\alpha_0+\alpha_1+2\alpha_2+2\alpha_3+2\alpha_4+\alpha_5+\alpha_6)(t-1)^2}, \quad K_2:=-H_2.
\end{equation}
Two Hamiltonians $K_1$ and $K_2$ satisfy
\begin{equation}
\{K_1,K_2\}+\left(\frac{\partial}{\partial s} \right)K_1-\left(\frac{\partial}{\partial t} \right)K_2=0,
\end{equation}
where  $\{,\}$ denotes the Poisson brackets:
\begin{equation}
\{L_1,L_2\}=\frac{\partial L_1}{\partial p_1}\frac{\partial L_2}{\partial q_1}-\frac{\partial L_1}{\partial q_1}\frac{\partial L_2}{\partial p_1}+\frac{\partial L_1}{\partial p_2}\frac{\partial L_2}{\partial q_2}-\frac{\partial L_1}{\partial q_2}\frac{\partial L_2}{\partial p_2}.
\end{equation}
\end{proposition}

We note that the system \eqref{eq:D6} admits affine Weyl group symmetry of type $D_6^{(1)}$ as the group of its B{\"a}cklund transformations, whose generators $s_0,s_1,\ldots,s_6,{\pi}_j \ (j=1,2,3)$ defined as follows$:$ with {\it the notation} $(*):=(q_1,p_1,q_2,p_2,t,s;\alpha_0,\alpha_1,\ldots,\alpha_6)$\rm{: \rm}
{\small
\begin{align}\label{symmetry:D6}
\begin{split}
s_0:(*) \rightarrow &\left(q_1,p_1-\frac{\alpha_0}{q_1-t},q_2,p_2,t,s;-\alpha_0,\alpha_1,\alpha_2+\alpha_0,\alpha_3,\alpha_4,\alpha_5,\alpha_6 \right),\\
s_1:(*) \rightarrow &\left(q_1,p_1-\frac{\alpha_1}{q_1-s},q_2,p_2,t,s;\alpha_0,-\alpha_1,\alpha_2+\alpha_1,\alpha_3,\alpha_4,\alpha_5,\alpha_6 \right),\\
s_2:(*) \rightarrow &\left(q_1+\frac{\alpha_2}{p_1},p_1,q_2,p_2,t,s;\alpha_0+\alpha_2,\alpha_1+\alpha_2,-\alpha_2,\alpha_3+\alpha_2,\alpha_4,\alpha_5,\alpha_6 \right),\\
s_3:(*) \rightarrow &\left(q_1,p_1-\frac{\alpha_3}{q_1-q_2},q_2,p_2+\frac{\alpha_3}{q_1-q_2},t,s;\alpha_0,\alpha_1,\alpha_2+\alpha_3,-\alpha_3,\alpha_4+\alpha_3,\alpha_5,\alpha_6 \right),\\
s_4:(*) \rightarrow &\left(q_1,p_1,q_2+\frac{\alpha_4}{p_2},p_2,t,s;;\alpha_0,\alpha_1,\alpha_2,\alpha_3+\alpha_4,-\alpha_4,\alpha_5+\alpha_4,\alpha_6+\alpha_4 \right),\\
s_5:(*) \rightarrow &\left(q_1,p_1,q_2,p_2-\frac{\alpha_5}{q_2-1},t,s;\alpha_0,\alpha_1,\alpha_2,\alpha_3,\alpha_4+\alpha_5,-\alpha_5,\alpha_6 \right),\\
s_6:(*) \rightarrow &\left(q_1,p_1,q_2,p_2-\frac{\alpha_6}{q_2},t,s;\alpha_0,\alpha_1,\alpha_2,\alpha_3,\alpha_4+\alpha_6,\alpha_5,-\alpha_6 \right),\\
        {\pi}_1: (*) \rightarrow  &(\frac{(t-1)q_1}{t-q_1-t s+t s q_1},\frac{(-t+q_1+t s-t s q_1)(t p_1-q_1 p_1-\alpha_2-t s p_1+t s q_1 p_1+\alpha_2 t s)}{t(t-1)(s-1)},\\
&\frac{(t-1)q_2}{t-q_2-t s+t s q_2},\frac{(-t+q_2+t s-t s q_2)(t p_2-q_2 p_2-\alpha_4-t s p_2+t s q_2 p_2+\alpha_4 t s)}{t(t-1)(s-1)},\\
&\frac{s(t-1)}{t-s-t s+t s^2},\frac{1}{s};\alpha_1,\alpha_0,\alpha_2,\alpha_3,\alpha_4,\alpha_5,\alpha_6 ), \\
        {\pi}_2: (*) \rightarrow  &(1-q_1,-p_1,1-q_2,-p_2,1-t,1-s;\alpha_0,\alpha_1,\alpha_2,\alpha_3,\alpha_4,\alpha_6,\alpha_5), \\
        {\pi}_3: (*) \rightarrow  &(\frac{t(q_2-s)}{t(q_2-s)+s^2 (t-q_2)},\frac{(t(q_2-s)+s^2 (t-q_2))(t(q_2-s)p_2+\alpha_4(t-s^2)+s^2 (t-q_2)p_2)}{t s^2 (t-s)},\\
&\frac{t(q_1-s)}{t(q_1-s)+s^2 (t-q_1)},\frac{(t(q_1-s)+s^2 (t-q_1))(t(q_1-s)p_1+\alpha_2(t-s^2)+s^2 (t-q_1)p_1)}{t s^2 (t-s)},\\
&-\frac{(s-1)t}{t-t s+s^2(t-1)},-\frac{1}{s-1};\alpha_5,\alpha_6,\alpha_4,\alpha_3,\alpha_2,\alpha_0,\alpha_1 ).
\end{split}
\end{align}
}
We remark that the system \eqref{eq:D6} has the following invariant divisors\rm{:\rm}
\begin{center}
\begin{tabular}{|c|c|c|} \hline
parameter's relation & invariant divisor \\ \hline
$\alpha_0=0$ & $f_0:=q_1-t$  \\ \hline
$\alpha_1=0$ & $f_1:=q_1-s$  \\ \hline
$\alpha_2=0$ & $f_2:=p_1$  \\ \hline
$\alpha_3=0$ & $f_3:=q_1-q_2$  \\ \hline
$\alpha_4=0$ & $f_4:=p_2$  \\ \hline
$\alpha_5=0$ & $f_5:=q_2-1$  \\ \hline
$\alpha_6=0$ & $f_6:=q_2$  \\ \hline
\end{tabular}
\end{center}
We note that when $\alpha_2=0$, we see that the system \eqref{eq:D6} admits a particular solution $p_1=0$, and when $\alpha_3=0$, after we make the birational and symplectic transformation:
\begin{equation}
x_3=q_1-q_2, \ y_3=p_1, \ z_3=q_2, \ w_3=p_2+p_1
\end{equation}
we see that the system \eqref{eq:D6} admits a particular solution $x_3=0$.

The B{\"a}cklund transformations of the system of type $D_6^{(1)}$ satisfy
\begin{equation}
s_i(g)=g+\frac{\alpha_i}{f_i}\{f_i,g\}+\frac{1}{2!} \left(\frac{\alpha_i}{f_i} \right)^2 \{f_i,\{f_i,g\} \}+\cdots \quad (g \in {\mathbb C}(t)[q_1,p_1,q_2,p_2]),
\end{equation}
where $\{,\}$ is the Poisson bracket such that $\{p_i,q_j\}={\delta}_{ij}$ (see \cite{N2}).

Since these B{\"a}cklund transformations have Lie theoretic origin, similarity reduction of a Drinfeld-Sokolov hierarchy admits such a B{\"a}cklund symmetry (see \cite{Yamada}).

\section{Appendix A: Accessible singularities of the system \eqref{eq:D6}}

\begin{figure}
\unitlength 0.1in
\begin{picture}( 34.9400, 18.6000)( 16.7000,-22.7000)
%
\special{pn 8}%
\special{pa 2430 888}%
\special{pa 1758 2190}%
\special{fp}%
\special{pa 1758 2190}%
\special{pa 5164 2190}%
\special{fp}%
%
\special{pn 8}%
\special{pa 2424 888}%
\special{pa 4432 888}%
\special{fp}%
%
\special{pn 8}%
\special{pa 4418 888}%
\special{pa 5158 2190}%
\special{fp}%
%
\special{pn 8}%
\special{pa 2430 888}%
\special{pa 2852 1326}%
\special{fp}%
\special{pa 2852 1326}%
\special{pa 2604 1842}%
\special{fp}%
%
\special{pn 8}%
\special{pa 1758 2182}%
\special{pa 2604 1834}%
\special{fp}%
%
\special{pn 8}%
\special{pa 2844 1326}%
\special{pa 4106 1326}%
\special{fp}%
\special{pa 4106 1326}%
\special{pa 4418 888}%
\special{fp}%
%
\special{pn 8}%
\special{pa 4112 1318}%
\special{pa 4344 1850}%
\special{fp}%
\special{pa 4344 1850}%
\special{pa 5152 2190}%
\special{fp}%
%
\special{pn 8}%
\special{pa 2610 1826}%
\special{pa 4324 1826}%
\special{fp}%
\put(16.7000,-24.4000){\makebox(0,0)[lb]{$Y_3$}}%
\put(50.5200,-24.2500){\makebox(0,0)[lb]{$W_4$}}%
\put(30.1100,-21.0100){\makebox(0,0)[lb]{$C_3$}}%
%
\special{pn 20}%
\special{pa 1924 1890}%
\special{pa 2712 1616}%
\special{fp}%
%
\special{pn 20}%
\special{pa 4752 2012}%
\special{pa 4258 1114}%
\special{fp}%
\put(21.1600,-17.6200){\makebox(0,0)[lb]{$C_0$}}%
\put(45.0500,-15.7500){\makebox(0,0)[lb]{$C_5$}}%
%
\special{pn 8}%
\special{pa 2430 888}%
\special{pa 2858 1332}%
\special{fp}%
%
\special{pn 20}%
\special{pa 4418 888}%
\special{pa 5158 2198}%
\special{fp}%
\put(47.0400,-13.9700){\makebox(0,0)[lb]{$C_6$}}%
\put(22.5000,-14.3000){\makebox(0,0)[lb]{$C_1$}}%
%
\special{pn 8}%
\special{pa 2438 896}%
\special{pa 3504 410}%
\special{dt 0.045}%
%
\special{pn 8}%
\special{pa 3518 434}%
\special{pa 4444 888}%
\special{dt 0.045}%
%
\special{pn 8}%
\special{pa 3490 418}%
\special{pa 3490 1502}%
\special{dt 0.045}%
%
\special{pn 8}%
\special{pa 2618 1818}%
\special{pa 3498 1510}%
\special{dt 0.045}%
%
\special{pn 8}%
\special{pa 3490 1510}%
\special{pa 4344 1826}%
\special{dt 0.045}%
%
\special{pn 8}%
\special{pa 2870 1326}%
\special{pa 3498 1074}%
\special{dt 0.045}%
%
\special{pn 8}%
\special{pa 3490 1082}%
\special{pa 4124 1326}%
\special{dt 0.045}%
%
\special{pn 20}%
\special{pa 3172 2182}%
\special{pa 3458 1326}%
\special{fp}%
\put(16.8000,-19.6000){\makebox(0,0)[lb]{$r_0$}}%
\put(18.2000,-14.7000){\makebox(0,0)[lb]{$r_1$}}%
\put(29.9000,-23.8000){\makebox(0,0)[lb]{$r_3$}}%
\put(46.6000,-18.6000){\makebox(0,0)[lb]{$r_5$}}%
\put(50.7000,-20.6000){\makebox(0,0)[lb]{$r_6$}}%
\put(32.7000,-8.2000){\makebox(0,0)[lb]{$r_2$}}%
\put(35.2000,-14.7000){\makebox(0,0)[lb]{$r_4$}}%
%
\special{pn 20}%
\special{pa 2160 1450}%
\special{pa 2790 1450}%
\special{fp}%
\end{picture}%
\label{CPD6fig4}
\caption{This figure denotes the boundary divisor ${\mathcal H}$ of ${\mathcal S}$ (see \cite{sasa3}). The bold lines $C_i \ i=0,1,3,5,6$ (see \eqref{accessible sing.}) in ${\mathcal H}$ denote the accessible singular loci of the system \eqref{eq:D6}. }
\end{figure}
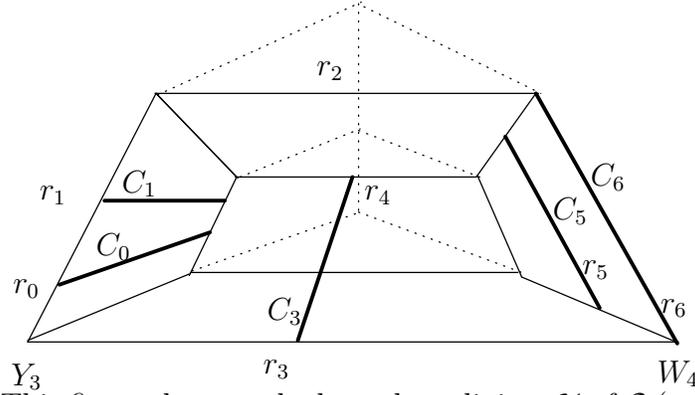

In order to consider the singularity analysis for the system \eqref{eq:D6}, as a compactification of ${\mathbb C}^4$ which is the phase space of the system \eqref{eq:D6}, we take 4-dimensional complex manifold $\mathcal S$ given in the paper \cite{sasa3}. This manifold can be considered as a generalization of the Hirzebruch surface.

We easily see that the rational vector field $\tilde v$ associated with the system \eqref{eq:D6} satisfies the condition:
\begin{equation}
\tilde v \in H^0({\mathcal S},\Theta_{\mathcal S}(-\log{\mathcal H})({\mathcal H})).
\end{equation}
We remark that in each coordinate system $(X_i,Y_i,Z_i,W_i), \ (i=1,2,5)$ (see \eqref{patchcoordinate}), the right hand-sides of the system \eqref{eq:D6} are polynomial.

We also remark that this rational vector field associated with  the system \eqref{eq:D6} has a pole along only divisor ${\mathcal H}$, whose order is one.

The rational vector field $\tilde v$ associated with the system \eqref{eq:D6}  has five accessible singular curves $C_i \cong {\mathbb P}^1 \ (i=0,1,3,5,6)$:
\begin{equation}\label{accessible sing.}
  \left\{
  \begin{aligned}
    C_0=&\{(X_3,Y_3,Z_3,W_3)|X_3=t,Y_3=W_3=0\}\\
& \cup \left\{(X_8,Y_8,Z_8,W_8)|X_8=\frac{1}{t},Y_8=W_8=0 \right\} \cong {\mathbb P}^1,\\
    C_1=&\{(X_3,Y_3,Z_3,W_3)|X_3=s,Y_3=W_3=0\}\\
& \cup \left\{(X_8,Y_8,Z_8,W_8)|X_8=\frac{1}{s},Y_8=W_8=0 \right\} \cong {\mathbb P}^1,\\
    C_3=&\{(X_3,Y_3,Z_3,W_3)|X_3=Z_3,Y_3=0,W_3=-1\}\\
& \cup \{(X_8,Y_8,Z_8,W_8)|X_8=Z_8,Y_8=0,W_8=-1\} \cong {\mathbb P}^1,\\
    C_5=&\{(X_4,Y_4,Z_4,W_4)|Y_4=0,Z_4=1,W_4=0\}\\
& \cup \{(X_9,Y_9,Z_9,W_9)|Y_9=0,Z_9=1,W_9=0\} \cong {\mathbb P}^1,\\
    C_6=&\{(X_4,Y_4,Z_4,W_4)|Y_4=Z_4=W_4=0\}\\
& \cup \{(X_7,Y_7,Z_7,W_7)|Y_7=Z_7=W_7=0\} \cong {\mathbb P}^1,
   \end{aligned}
  \right. 
\end{equation}

\begin{tabular}{|c||c|c|} \hline 
Singular curve & $C_0$ & $C_1$ \\ \hline
Ratio of local index    &  [2,0,1] (dim. 3) & [2,0,1] (dim. 3)   \\ \hline
\end{tabular}

\begin{tabular}{|c||c|c|c|} \hline 
Singular curve & $C_3$ & $C_5$ & $C_6$ \\ \hline
Ratio of local index & [2,0,1] (dim. 3) &  [0,1,2] (dim. 3)  & [0,1,2] (dim. 3)   \\ \hline
\end{tabular}

Here, the coordinate systems $(X_i,Y_i,Z_i,W_i) \ (i=1,2,\ldots,11)$ (see  \cite{sasa3}) are given by
{\small
\begin{align}\label{patchcoordinate}
\begin{split}
&(X_1,Y_1,Z_1,W_1)=\left(\frac{1}{q_1},-(p_1q_1+\alpha_2)q_1,q_2,p_2 \right), \ (X_2,Y_2,Z_2,W_2)=\left(q_1,p_1,\frac{1}{q_2},-(p_2q_2+\alpha_4)q_2 \right),\\
&(X_3,Y_3,Z_3,W_3)=\left(q_1,\frac{1}{p_1},q_2,\frac{p_2}{p_1} \right), \ (X_4,Y_4,Z_4,W_4)=\left(q_1,\frac{p_1}{p_2},q_2,\frac{1}{p_2} \right),\\
&(X_5,Y_5,Z_5,W_5)=\left(\frac{1}{q_1},-(p_1q_1+\alpha_2)q_1,\frac{1}{q_2},-(p_2q_2+\alpha_4)q_2 \right),\\
&(X_6,Y_6,Z_6,W_6)=\left(\frac{1}{q_1},-\frac{1}{(q_1p_1+\alpha_2)q_1},q_2,-\frac{p_2}{(q_1p_1+\alpha_2)q_1} \right),\\
&(X_7,Y_7,Z_7,W_7)=\left(\frac{1}{q_1},-\frac{(p_1q_1+\alpha_2)q_1}{p_2},q_2,\frac{1}{p_2} \right),\\
&(X_8,Y_8,Z_8,W_8)=\left(\frac{1}{q_1},-\frac{1}{(p_1q_1+\alpha_2)q_1},\frac{1}{q_2},\frac{(p_2q_2+\alpha_4)q_2}{(p_1q_1+\alpha_2)q_1} \right),\\
&(X_9,Y_9,Z_9,W_9)=\left(\frac{1}{q_1},\frac{(p_1q_1+\alpha_2)q_1}{(p_2q_2+\alpha_4)q_2},\frac{1}{q_2},-\frac{1}{(p_2q_2+\alpha_4)q_2} \right),\\
&(X_{10},Y_{10},Z_{10},W_{10})=\left(q_1,\frac{1}{p_1},\frac{1}{q_2},-\frac{(p_2q_2+\alpha_4)q_2}{p_1} \right),\\
&(X_{11},Y_{11},Z_{11},W_{11})=\left(q_1,-\frac{p_1}{(p_2q_2+\alpha_4)q_2},\frac{1}{q_2},-\frac{1}{(p_2q_2+\alpha_4)q_2} \right).
\end{split}
\end{align}
}

\section{Appendix B: Painlev\'e scheme of $D_6^{(1)}$ system}

Let us consider the system of the first order ordinary differential equations of polynomial type:
\begin{equation}\label{sys:1}
  \left\{
  \begin{aligned}
   \frac{\partial x}{\partial t} &=f_1(x,y,z,w), \quad \frac{\partial y}{\partial t} =f_2(x,y,z,w),\\
   \frac{\partial z}{\partial t} &=f_3(x,y,z,w), \quad \frac{\partial w}{\partial t} =f_4(x,y,z,w),
    \end{aligned}
  \right. 
\end{equation}
where $f_i \in {\mathbb C}(t,s)[x,y,z,w]$ and $deg(f_i)=5$ with respect to $x,y,z,w$.

We assume that associated vector field $v$
$$
v=\frac{\partial}{\partial t}+f_1(x,y,z,w)\frac{\partial}{\partial x}+f_2(x,y,z,w)\frac{\partial}{\partial y}+f_3(x,y,z,w)\frac{\partial}{\partial z}+f_4(x,y,z,w)\frac{\partial}{\partial w}
$$
belongs in 
\begin{equation}\label{conv}
v \in H^0({\mathcal S},\Theta_{\mathcal S}(-\log{\mathcal H})({\mathcal H})).
\end{equation}
Here, the complex manifold $\mathcal S$ is given in the paper \cite{sasa3} (see Section 4). This manifold can be considered as a generalization of the Hirzebruch surface.

This condition is equivalent to the following:
\begin{enumerate}\label{2}
\item Holomorphy in the coordinate systems $(x_1,y_1,z_1,w_1)=(1/x,-(xy+\alpha_2)x,z,w),\\
(x_2,y_2,z_2,w_2)=(x,y,1/z,-(wz+\alpha_4)z)$.
\item In the coordinate system $(X,Y,Z,W)=(x,1/y,z,w/y)$, the differential system \eqref{sys:1} must be taken of the form:
\begin{equation}
  \left\{
  \begin{aligned}
   \frac{\partial X}{\partial t} &=\frac{F_1(X,Y,Z,W)}{Y},  \ \frac{\partial Y}{\partial t} =F_2(X,Y,Z,W),\\
   \frac{\partial Z}{\partial t} &=\frac{F_3(X,Y,Z,W)}{Y}, \ \frac{\partial W}{\partial t} =\frac{F_4(X,Y,Z,W)}{Y} \quad (F_i \in {\mathbb C}(t,s)[X,Y,Z,W]).
   \end{aligned}
  \right. 
\end{equation}
\item In the coordinate system $(X,Y,Z,W)=(x,y/w,z,1/w)$, the differential system \eqref{sys:1} must be taken of the form:
\begin{equation}
  \left\{
  \begin{aligned}
   \frac{\partial X}{\partial t} &=\frac{G_1(X,Y,Z,W)}{W},  \ \frac{\partial Y}{\partial t} =\frac{G_2(X,Y,Z,W)}{W},\\
   \frac{\partial Z}{\partial t} &=\frac{G_3(X,Y,Z,W)}{W}, \ \frac{\partial W}{\partial t} =G_4(X,Y,Z,W) \quad (G_i \in {\mathbb C}(t,s)[X,Y,Z,W]).
   \end{aligned}
  \right. 
\end{equation}
\end{enumerate}

We easily see that the system \eqref{sys:1} satisfying the condition \eqref{conv} has forty undetermined coefficients.

For the system \eqref{sys:1} satisfying the condition \eqref{conv}, we give the following Painlev\'e scheme (see \cite{sasa8}):

{\small
\begin{equation}\label{asd}
\begin{pmatrix}
{\rm Accessible \ singular \ loci} & C_0 & C_1 & C_3 &  C_5 &  C_6 \\
{\rm (In \ each \ coordinate \ origin} &  P_0 \in C_0 & P_1 \in C_1 & P_3 \in C_3 &  P_5 \in C_5 &  P_6 \in C_6, )\\
{\rm Continued \  Ratio \ of \ Local \ index} & \begin{bmatrix}
n_1\\
1 \\
0
\end{bmatrix} & \begin{bmatrix}
n_2\\
1 \\
0
\end{bmatrix} & \begin{bmatrix}
n_5\\
n_6 \\
0
\end{bmatrix} &\begin{bmatrix}
n_4\\
1 \\
0
\end{bmatrix} & \begin{bmatrix}
n_3\\
1 \\
0
\end{bmatrix} 
\end{pmatrix} \ (n_i \in {\mathbb C}),
\end{equation} }
where each of the curves $C_0,C_1,C_5,C_6$ is a ${\mathbb P}^1$-fiber in the surface ${\mathbb P}^1 \times {\mathbb P}^1$ and the curve $C_3$ is 
its diagonal component.

Here, \begin{equation*}
  \left\{
  \begin{aligned}
    C_0=&\{(X_3,Y_3,Z_3,W_3)|X_3=t,Y_3=W_3=0\} \subset {\mathbb P}^1,\\
    C_1=&\{(X_3,Y_3,Z_3,W_3)|X_3=s,Y_3=W_3=0\} \subset {\mathbb P}^1,\\
    C_3=&\{(X_3,Y_3,Z_3,W_3)|X_3=Z_3,Y_3=0,W_3=-1\} \subset {\mathbb P}^1,\\
    C_5=&\{(X_4,Y_4,Z_4,W_4)|Y_4=0,Z_4=1,W_4=0\} \subset {\mathbb P}^1,\\
    C_6=&\{(X_4,Y_4,Z_4,W_4)|Y_4=Z_4=W_4=0\} \subset {\mathbb P}^1,
   \end{aligned}
  \right. 
\end{equation*} and \begin{equation*}
  \left\{
  \begin{aligned}
    P_0=&\{(X_3,Y_3,Z_3,W_3)=(t,0,0,0)\},\\
    P_1=&\{(X_3,Y_3,Z_3,W_3)=(s,0,0,0)\},\\
    P_3=&\{(X_3,Y_3,Z_3,W_3)=(0,0,0,-1)\},\\
    P_5=&\{(X_4,Y_4,Z_4,W_4)=(0,0,1,0)\},\\
    P_6=&\{(X_4,Y_4,Z_4,W_4)=(0,0,0,0)\}.
   \end{aligned}
  \right. 
\end{equation*}
The coordinate systems $(X_i,Y_i,Z_i,W_i) \ (i=1,2,\ldots,11)$ (cf.  \cite{sasa3}) are given in \eqref{patchcoordinate} in Section 4.

For example, the eigenvalues of the matrix of linear approximation around the accessible singular point $P_0$ is given by $(n_1 f(t,s),f(t,s), 0, f(t,s))$, whose continued  ratio is
\begin{equation}
\left(\frac{n_1 f(t,s)}{f(t,s)},\frac{0}{f(t,s)},\frac{f(t,s)}{f(t,s)}, \right)=(n_1,0,1),
\end{equation}
where $f(t,s) \in {\mathbb C}(t,s)$.

Then, we can obtain the following {\bf eigenvalues's relations.}
\begin{proposition}
The eigenvalues $n_i, \ (i=1,2,\ldots,6)$ in \eqref{asd} satisfy the following relations\rm{: \rm}
\begin{equation}\label{eigenvalues}
  \left\{
  \begin{aligned}
   &n_1 n_2(2n_3 n_4-n_3-n_4)-n_5(2n_1 n_2-n_1-n_2)(n_3-1)n_4=0,\\
   &(n_3-1)n_6=1,
   \end{aligned}
  \right. 
\end{equation}
\end{proposition}
where $D_6^{(1)}$ case is $(n_1,n_2,n_3,n_4,n_5,n_6)=(2,2,2,2,2,1)$.

We remark that by making a change of variables
\begin{equation}
N_5=\frac{2n_1 n_2 n_3 n_4 n_5}{2n_1 n_2 n_3(n_4-n_5)-(n_1+n_2)n_3 n_4 n_5-n_1 n_2 n_3-n_1 n_2 n_4-n_1 n_4 n_5-n_2 n_4 n_5}
\end{equation}
the first equation in \eqref{eigenvalues} can be transformed into the equation (cf. \cite{cos1,cos2,sasa8})
\begin{equation}
\frac{1}{n_1}+\frac{1}{n_2}+\frac{1}{n_3}+\frac{1}{n_4}+\frac{1}{N_5}=1.
\end{equation}
Here, the eigenvalue's relation;
\begin{equation}
\frac{1}{n_1}+\frac{1}{n_2}+\frac{1}{n_3}+\frac{1}{n_4}+\frac{1}{n_5}=b \quad (b \in {\mathbb C})
\end{equation}
can be transformed into a one-parameter family of quintic hypersurfaces (see \cite{Candelas});
\begin{equation}
x_1^5+x_2^5+x_3^5+x_4^5+x_5^5-b x_1 x_2 x_3 x_4 x_5=0,
\end{equation}
where we can make a change of variables:
\begin{align}
\begin{split}
&n_1=\frac{x_2 x_3 x_4 x_5}{x_1^4}, \quad n_2=\frac{x_1 x_3 x_4 x_5}{x_2^4}, \quad n_3=\frac{x_1 x_2 x_4 x_5}{x_3^4}, \quad n_4=-\frac{x_1 x_2 x_3 x_4 x_5 n_5}{x_1 x_2 x_3 x_4 x_5-(x_4^5+x_5^5) n_5}.
\end{split}
\end{align}

In general, the eigenvalue's relation (cf. \cite{cos1,cos2});
\begin{equation}
\frac{1}{n_1}+\frac{1}{n_2}+\cdots+\frac{1}{n_N}=b \quad (b \in {\mathbb C})
\end{equation}
can be transformed into a one-parameter family of Calabi-Yau hypersurfaces;
\begin{equation}
x_1^N+x_2^N+\cdots+x_N^N-b x_1 x_2 \cdots x_N=0,
\end{equation}
where we can make a change of variables:
\begin{align}
\begin{split}
&n_1=\frac{x_2 x_3 \cdots x_N}{x_1^{N-1}}, \quad n_2=\frac{x_1 x_3 \cdots x_N}{x_2^{N-1}}, \cdots ,n_{N-2}=\frac{x_1 \cdots x_{N-3} x_{N-1} x_N}{x_{N-2}^{N-1}},\\
&n_{N-1}=-\frac{x_1 x_2 \cdots x_N n_N}{x_1 x_2 \cdots x_N-(x_{N-1}^N+x_N^N) n_N}.
\end{split}
\end{align}
We will call these eigenvalues's relations {\it Painlev\'e-Cosgrove relations} (cf. \cite{cos1,cos2,sasa8}).

\section{Appendix C: Holomorphy of the system of type $D_6^{(1)}$}

In this appendix, we list some holomorphy conditions of the system of type $D_6^{(1)}$.

{\bf Hamiltonians $H_{0}^{(1)}=r_{0}(H_1^{D_6^{(1)}}-p_1), \ H_{0}^{(2)}=r_{0}(H_2^{D_6^{(1)}}), \ r_{0}:x=-((q_1-t)p_1-\alpha_0)p_1, \ y=\frac{1}{p_1}, \ z=q_2, \ w=p_2$}

\begin{align*}
&r_0^{0}:x_0=-(q_1p_1-\alpha_0)p_1, \ y_0=\frac{1}{p_1}, \ z_0=q_2, \ w_0=p_2, \\
&r_1^{0}:x_1=q_1+\frac{\alpha_1-\alpha_0}{p_1}+\frac{s-t}{p_1^2}, \ y_1=p_1, \ z_1=q_2, \ w_1=p_2, \\
&r_2^{0}:x_2=-(q_1p_1-\alpha_2-\alpha_0)p_1, \ y_2=\frac{1}{p_1}, \ z_2=q_2 \ w_2=p_2, \\
&r_3^{0}:x_3=q_1+\frac{\alpha_3-\alpha_0}{p_1}+\frac{q_2-t}{p_1^2}, \ y_3=p_1, \ z_3=q_2, \ w_3=p_2+\frac{1}{p_1}, \\
&r_4^{0}:x_4=q_1, \ y_4=p_1, \ z_4=\frac{1}{q_2}, \ w_4=-(p_2q_2+\alpha_4)q_2, \\
&r_5^{0}:x_5=q_1, \ y_5=p_1, \ z_5=-((q_2-1)p_2-\alpha_5)p_2, \ w_5=\frac{1}{p_2}, \\
&r_6^{0}:x_6=q_1, \ y_6=p_1, \ z_6=-(q_2p_2-\alpha_6)p_2, \ w_6=\frac{1}{p_2},
\end{align*}
where $r_3^{0} \left(H_{0}^{(1)}+\frac{1}{p_1} \right), \ r_1^{0} \left(H_{0}^{(1)}+\frac{1}{p_1} \right), \ r_1^{0} \left(H_{0}^{(2)}-\frac{1}{p_1} \right)$. Here, for notational convenience, we have renamed $(x,y,z,w)$ to $(q_1,p_1,q_2,p_2)$ (which are not the same as the previous $(q_1,p_1,q_2,p_2)$). It is still an open question whether the transformation $r_3^{0}$ can be considered as its auto-B{\"a}cklund transformation or not.

We remark that each of transformations like $r_3^{0}$ is a B{\"a}cklund transformation for the Noumi-Yamada system of type $A_5^{(1)}$ and Garnier system in two variables, respectively.

We note that we will construct the transformation $r_3^{0}$.

{\bf Step 1:} We make a change of variables:
$$
X^{(1)}=q_1 p_1^2, \quad Y^{(1)}=p_1, \quad Z^{(1)}=q_2, \quad W^{(1)}=p_2 p_1.
$$

{\bf Step 2:} We blow up along the curve $\{(X^{(1)},Y^{(1)},Z^{(1)},W^{(1)})|X^{(1)}=t-Z^{(1)}, \ Y^{(1)}=0, \  W^{(1)}=-1\}$:
$$
X^{(2)}=\frac{X^{(1)}-(t-Z^{(1)})}{Y^{(1)}}, \quad Y^{(2)}=Y^{(1)}, \quad Z^{(2)}=Z^{(1)}, \quad W^{(2)}=\frac{W^{(1)}+1}{Y^{(1)}}.
$$

{\bf Step 3:} We blow up along the surface $\{(X^{(2)},Y^{(2)},Z^{(2)},W^{(2)})|X^{(2)}=-(\alpha_3-\alpha_0), \ Y^{(2)}=0 \}$:
$$
X^{(3)}=\frac{X^{(2)}+\alpha_3-\alpha_0}{Y^{(2)}}, \quad Y^{(3)}=Y^{(2)}, \quad Z^{(3)}=Z^{(2)}, \quad W^{(3)}=W^{(2)}.
$$
By taking the coordinate system as
$$
(x_3,y_3,z_3,w_3)=(X^{(3)},Y^{(3)},Z^{(3)},W^{(3)}),
$$
we can obtain the coordinate system $r_3^{0}$.

\vspace{0.5cm}

{\bf Hamiltonians $H_{6}^{(1)}=r_{6}(H_1^{D_6^{(1)}}), \ H_{6}^{(2)}=r_{6}(H_2^{D_6^{(1)}}), \ r_{6}:x=q_1, \ y=p_1, \ z=-(q_2p_2-\alpha_6)p_2, \ w=\frac{1}{p_2}$}

{\footnotesize
\begin{align*}
&r_0^{6}:x_0=-((q_1-t)p_1-\alpha_0)p_1, \ y_0=\frac{1}{p_1}, \ z_0=q_2, \ w_0=p_2, \\
&r_1^{6}:x_1=-((q_1-s)p_1-\alpha_1)p_1, \ y_1=\frac{1}{p_1}, \ z_1=q_2, \ w_1=p_2, \\
&r_2^{6}:x_2=\frac{1}{q_1}, \ y_2=-(p_1q_1+\alpha_2)q_1, \ z_2=q_2 \ w_2=p_2, \\
&r_3^{6}:x_3=q_1, \ y_3=p_1+\frac{1}{p_2}, \ z_3=q_2+\frac{\alpha_3-\alpha_6}{p_2}+\frac{q_1}{p_2^2}, \ w_3=p_2, \\
&r_4^{6}:x_4=q_1, \ y_4=p_1, \ z_4=-(q_2p_2-\alpha_4-\alpha_6)p_2, \ w_4=\frac{1}{p_2}, \\
&r_5^{6}:x_5=q_1, \ y_5=p_1, \ z_5=q_2+\frac{\alpha_5-\alpha_6}{p_2}+\frac{1}{p_2^2}, \ w_5=p_2, \\
&r_6^{6}:x_6=q_1, \ y_6=p_1, \ z_6=-(q_2p_2-\alpha_6)p_2, \ w_6=\frac{1}{p_2},
\end{align*}
where $r_0^{6} \left(H_{6}^{(1)}-p_1 \right), \ r_1^{6} \left(H_{6}^{(2)}-p_1 \right)$.
}

We remark that by making a change of variables $(q_i,p_i)$ and $\alpha_j$, the following transformation $s_5^{6}$ associated with $r_5^{6}$ becomes a B{\"a}cklund transformation:
\begin{align}
\begin{split}
        s_5^{6}: (*) &\rightarrow \left(-q_1,-p_1,-\left(q_2+\frac{\alpha_5-\alpha_6}{p_2}+\frac{1}{p_2^2} \right),-p_2,1-t,1-s;\alpha_0,\alpha_1,\alpha_2,\alpha_3,\alpha_4,\alpha_6,\alpha_5 \right).
        \end{split}
        \end{align}

We remark that for the Hamiltonian system with polynomial Hamiltonians $H_{6}^{(1)},H_{6}^{(2)}$ we can obtain another holomorphy conditions explicitly given by $r_i^{6} \ (i=0,1,2,4,5,6)$ and $\tilde{r}_3^{6}:x_3=-\left( \left(q_1+\frac{q_2}{p_1^2} \right)p_1-\frac{2q_2(p_2p_1+1)}{p_1}-(\alpha_3-\alpha_6) \right)p_1, \ y_3=\frac{1}{p_1},z_3=\frac{q_2}{p_1^2}, \ w_3=(p_2p_1+1)p_1$.

By these conditions, we can recover the polynomial Hamiltonians $H_{6}^{(1)},H_{6}^{(2)}$.

We note that we will construct the transformation $\tilde{r}_3^{6}$.

{\bf Step 1:} We make a change of variables:
$$
X^{(1)}=q_1, \quad Y^{(1)}=\frac{1}{p_1}, \quad Z^{(1)}=\frac{q_2}{p_1^2}, \quad W^{(1)}=p_2 p_1.
$$

{\bf Step 2:} We blow up along the curve $\{(X^{(1)},Y^{(1)},Z^{(1)},W^{(1)})|X^{(1)}=-Z^{(1)}, \ Y^{(1)}=0, \  W^{(1)}=-1\}$:
$$
X^{(2)}=\frac{X^{(1)}+Z^{(1)}}{Y^{(1)}}, \quad Y^{(2)}=Y^{(1)}, \quad Z^{(2)}=Z^{(1)}, \quad W^{(2)}=\frac{W^{(1)}+1}{Y^{(1)}}.
$$

{\bf Step 3:} We blow up along the surface $\{(X^{(2)},Y^{(2)},Z^{(2)},W^{(2)})|X^{(2)}=2Z^{(2)} W^{(2)}+\alpha_3-\alpha_6, \ Y^{(2)}=0 \}$:
$$
X^{(3)}=\frac{X^{(2)}-(2 Z^{(2)} W^{(2)}+\alpha_3-\alpha_6)}{Y^{(2)}}, \quad Y^{(3)}=Y^{(2)}, \quad Z^{(3)}=Z^{(2)}, \quad W^{(3)}=W^{(2)}.
$$
By choosing a new coordinate system as
$$
(x_3,y_3,z_3,w_3)=(-X^{(3)},Y^{(3)},Z^{(3)},W^{(3)}),
$$
we can obtain the coordinate system $\tilde{r}_3^{6}$.

\vspace{0.5cm}

{\bf Hamiltonians $H_{3}^{(1)}=\tilde{r}_{3}(H_1^{D_6^{(1)}}),\ H_{3}^{(2)}=\tilde{r}_{3}(H_2^{D_6^{(1)}}),\ \tilde{r}_{3}:x=q_1, \ y=p_1+p_2, \ z=-((q_2-q_1)p_2-\alpha_3)p_2, \ w=\frac{1}{p_2}$}

\begin{align*}
&r_0^{3}:x_0=-((q_1-t)p_1-\alpha_0)p_1, \ y_0=\frac{1}{p_1}, \ z_0=q_2, \ w_0=p_2, \\
&r_1^{3}:x_1=-((q_1-s)p_1-\alpha_1)p_1, \ y_1=\frac{1}{p_1}, \ z_1=q_2, \ w_1=p_2, \\
&r_2^{3}:x_2=q_1+\frac{\alpha_2 p_2}{p_1p_2-1}, \ y_2=p_1, \ z_2=q_2+\frac{\alpha_2 p_1}{p_1p_2-1} \ w_2=p_2, \\
&r_3^{3}:x_3=q_1, \ y_3=p_1, \ z_3=-(q_2p_2-\alpha_3)p_2, \ w_3=\frac{1}{p_2}, \\
&r_4^{3}:x_4=q_1, \ y_4=p_1, \ z_4=-(q_2p_2-\alpha_4-\alpha_3)p_2, \ w_4=\frac{1}{p_2}, \\
&r_5^{3}:x_5=q_1, \ y_5=p_1-\frac{1}{p_2}, \ z_5=q_2+\frac{\alpha_5-\alpha_3}{p_2}+\frac{1-q_1}{p_2^2}, \ w_5=p_2, \\
&r_6^{3}:x_6=q_1, \ y_6=p_1-\frac{1}{p_2}, \ z_6=q_2+\frac{\alpha_6-\alpha_3}{p_2}-\frac{q_1}{p_2^2}, \ w_6=p_2,
\end{align*}
where $r_0^{3} \left(H_{3}^{(1)}-p_1 \right), \ r_1^{3} \left(H_{3}^{(2)}-p_1 \right)$.

\vspace{0.5cm}

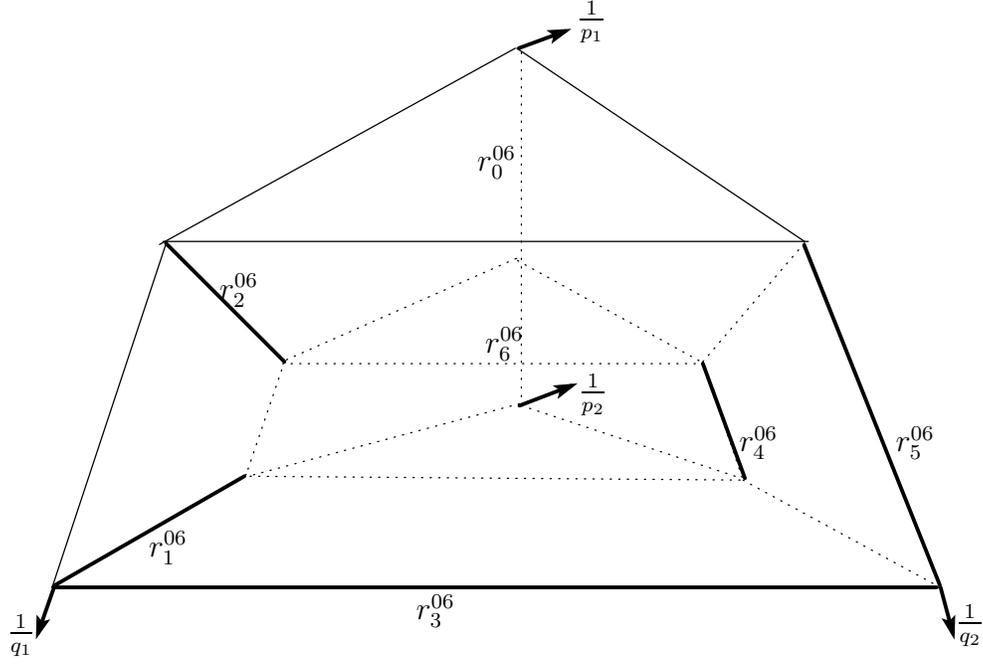
\begin{figure}
\unitlength 0.1in
\begin{picture}( 49.6700, 32.6400)(  7.6000,-34.3100)
%
\special{pn 8}%
\special{pa 3418 368}%
\special{pa 1556 1396}%
\special{fp}%
\special{pa 1574 1380}%
\special{pa 4952 1380}%
\special{fp}%
\special{pa 3426 368}%
\special{pa 4932 1380}%
\special{fp}%
%
\special{pn 8}%
\special{pa 1594 1380}%
\special{pa 994 3192}%
\special{fp}%
\special{pa 994 3192}%
\special{pa 5644 3192}%
\special{fp}%
\special{pa 4932 1388}%
\special{pa 5644 3192}%
\special{fp}%
%
\special{pn 20}%
\special{pa 1014 3192}%
\special{pa 5624 3192}%
\special{fp}%
%
\special{pn 20}%
\special{pa 3436 368}%
\special{pa 3680 278}%
\special{fp}%
\special{sh 1}%
\special{pa 3680 278}%
\special{pa 3610 282}%
\special{pa 3630 296}%
\special{pa 3624 320}%
\special{pa 3680 278}%
\special{fp}%
%
\special{pn 20}%
\special{pa 1004 3192}%
\special{pa 920 3432}%
\special{fp}%
\special{sh 1}%
\special{pa 920 3432}%
\special{pa 960 3376}%
\special{pa 938 3382}%
\special{pa 922 3362}%
\special{pa 920 3432}%
\special{fp}%
%
\special{pn 20}%
\special{pa 5644 3192}%
\special{pa 5708 3432}%
\special{fp}%
\special{sh 1}%
\special{pa 5708 3432}%
\special{pa 5710 3362}%
\special{pa 5694 3380}%
\special{pa 5672 3372}%
\special{pa 5708 3432}%
\special{fp}%
%
\special{pn 20}%
\special{pa 3440 2240}%
\special{pa 3712 2136}%
\special{fp}%
\special{sh 1}%
\special{pa 3712 2136}%
\special{pa 3642 2142}%
\special{pa 3662 2156}%
\special{pa 3656 2180}%
\special{pa 3712 2136}%
\special{fp}%
\put(37.4400,-3.3700){\makebox(0,0)[lb]{$\frac{1}{p_1}$}}%
\put(37.5000,-22.9000){\makebox(0,0)[lb]{$\frac{1}{p_2}$}}%
\put(57.2700,-34.9800){\makebox(0,0)[lb]{$\frac{1}{q_2}$}}%
\put(7.6000,-35.5000){\makebox(0,0)[lb]{$\frac{1}{q_1}$}}%
\put(29.0000,-34.0000){\makebox(0,0)[lb]{$r_3^{06}$}}%
%
\special{pn 8}%
\special{pa 1590 1400}%
\special{pa 2210 2010}%
\special{dt 0.045}%
\special{pa 2210 2010}%
\special{pa 2010 2600}%
\special{dt 0.045}%
\special{pa 2010 2600}%
\special{pa 990 3190}%
\special{dt 0.045}%
%
\special{pn 8}%
\special{pa 4930 1390}%
\special{pa 4400 2010}%
\special{dt 0.045}%
\special{pa 4400 2010}%
\special{pa 4610 2630}%
\special{dt 0.045}%
\special{pa 4610 2630}%
\special{pa 5620 3180}%
\special{dt 0.045}%
%
\special{pn 8}%
\special{pa 2200 2020}%
\special{pa 4380 2020}%
\special{dt 0.045}%
%
\special{pn 8}%
\special{pa 2030 2610}%
\special{pa 4620 2630}%
\special{dt 0.045}%
%
\special{pn 8}%
\special{pa 2220 2000}%
\special{pa 3420 1470}%
\special{dt 0.045}%
\special{pa 3420 1470}%
\special{pa 3420 1470}%
\special{dt 0.045}%
%
\special{pn 8}%
\special{pa 3430 1480}%
\special{pa 4390 2010}%
\special{dt 0.045}%
%
\special{pn 8}%
\special{pa 3450 390}%
\special{pa 3450 2230}%
\special{dt 0.045}%
%
\special{pn 8}%
\special{pa 3440 2230}%
\special{pa 2020 2610}%
\special{dt 0.045}%
%
\special{pn 8}%
\special{pa 3450 2240}%
\special{pa 4600 2620}%
\special{dt 0.045}%
%
\special{pn 20}%
\special{pa 1590 1390}%
\special{pa 2210 2010}%
\special{fp}%
%
\special{pn 20}%
\special{pa 4930 1400}%
\special{pa 5640 3180}%
\special{fp}%
\put(18.7000,-17.2000){\makebox(0,0)[lb]{$r_2^{06}$}}%
\put(45.9000,-25.2000){\makebox(0,0)[lb]{$r_4^{06}$}}%
\put(54.1000,-25.2000){\makebox(0,0)[lb]{$r_5^{06}$}}%
\put(32.2000,-10.5000){\makebox(0,0)[lb]{$r_0^{06}$}}%
\put(32.6000,-20.0000){\makebox(0,0)[lb]{$r_6^{06}$}}%
%
\special{pn 20}%
\special{pa 4400 2020}%
\special{pa 4620 2620}%
\special{fp}%
\put(15.0000,-30.7000){\makebox(0,0)[lb]{$r_1^{06}$}}%
%
\special{pn 20}%
\special{pa 1010 3180}%
\special{pa 2000 2610}%
\special{fp}%
\end{picture}%
\label{fig:D6r06'Ì"ÁˆÙ"_2}
\caption{Accessible singular loci of the system with Hamiltonians $H_{06}^{(1)},H_{06}^{(2)}$}
\end{figure}

{\bf Hamiltonians $H_{06}^{(1)}=r_{06}(H_1^{D_6^{(1)}}-p_1), \ H_{06}^{(2)}=r_{06}(H_2^{D_6^{(1)}}), \ r_{06}:x=-((q_1-t)p_1-\alpha_0)p_1, \ y=\frac{1}{p_1}, \ z=-(q_2p_2-\alpha_6)p_2, \ w=\frac{1}{p_2}$}
\begin{align*}
&r_0^{06}:x_0=-(q_1p_1-\alpha_0)p_1, \ y_0=\frac{1}{p_1}, \ z_0=q_2, \ w_0=p_2, \\
&r_1^{06}:x_1=q_1+\frac{\alpha_1-\alpha_0}{p_1}+\frac{s-t}{p_1^2}, \ y_1=p_1, \ z_1=q_2, \ w_1=p_2, \\
&r_2^{06}:x_2=-(q_1p_1-\alpha_2-\alpha_0)p_1, \ y_2=\frac{1}{p_1}, \ z_2=q_2 \ w_2=p_2, \\
&r_3^{06}:x_3=q_1+q_2+\frac{2q_2p_2+\alpha_3-\alpha_0-\alpha_6}{p_1}-\frac{t}{p_1^2}, \ y_3=p_1, \ z_3=q_2p_1^2, \ w_3=\frac{p_2+p_1}{p_1^2}, \\
&r_4^{06}:x_4=q_1, \ y_4=p_1, \ z_4=-(q_2p_2-\alpha_4-\alpha_6)p_2, \ w_4=\frac{1}{p_2}, \\
&r_5^{06}:x_5=q_1, \ y_5=p_1, \ z_5=q_2+\frac{\alpha_5-\alpha_6}{p_2}+\frac{1}{p_2^2}, \ w_5=p_2, \\
&r_6^{06}:x_6=q_1, \ y_6=p_1, \ z_6=-(q_2p_2-\alpha_6)p_2, \ w_6=\frac{1}{p_2},
\end{align*}
where $r_1^{06} \left(H_{06}^{(1)}+\frac{1}{p_1} \right), \ r_3^{06} \left(H_{06}^{(1)}+\frac{1}{p_1} \right), \ r_1^{06} \left(H_{06}^{(2)}-\frac{1}{p_1} \right)$.

We remark that the transformation $r_3^{06}$ is not its auto-B{\"a}cklund transformation. It is still an open question whether the transformation $r_3^{06}$ can be considered as the transformation denoted by the symbol $\odot$ in the Oshima's paper (see \cite{Oshima}).

It is also still an open question whether we can obtain the Hamiltonian system with $H_{06}^{(1)},H_{06}^{(2)}$ by solving $4 \times 4$ Lax pair  (cf. \cite{KNS,Oshima,Fuji}) satisfying the following Riemann scheme:
\begin{equation}
\begin{pmatrix}
X=0 & X=1 &  X=t &  X=\infty\\
\begin{matrix}
0 \\
0 \\
0 \\
\theta_1
\end{matrix}  & \begin{matrix}
0 \\
0 \\
\theta_2 \\
\theta_2
\end{matrix}  &  \begin{matrix}
0 \\
0 \\
\theta_3 \\
\theta_3
\end{matrix} &  \begin{matrix}
\alpha_0\\
\alpha_2+\alpha_0 \\
\alpha_6 \\
\alpha_4+\alpha_6
\end{matrix}
\end{pmatrix}.
\end{equation}

Here, we will conjecture the following relations between Riemann data and Holomorphy conditions $r_i^{06} \ (i=0,1,\ldots,6);$

{\footnotesize
$\begin{pmatrix}
X=0 \\
\begin{matrix}
0 \\
0 \\
0 \\
\theta
\end{matrix}
\end{pmatrix} \Longleftrightarrow $ Holomorphy condition $r_3^{06}$, \quad $\begin{pmatrix}
X=1,t \\
\begin{matrix}
0 \\
0 \\
\theta \\
\theta
\end{matrix}
\end{pmatrix} \Longleftrightarrow $ Holomorphy conditions $r_1^{06},r_5^{06}$,

$\begin{pmatrix}
X=\infty \\
\begin{matrix}
\alpha_0\\
\alpha_2+\alpha_0 \\
\alpha_6 \\
\alpha_4+\alpha_6
\end{matrix}
\end{pmatrix} \Longleftrightarrow $ Holomorphy conditions $\begin{pmatrix}
\begin{matrix}
r_0^{06} \\
r_2^{06} \\
r_6^{06} \\
r_4^{06}
\end{matrix}
\end{pmatrix}.$
}

We hope to call the Hamiltonians $H_{06}^{(1)},H_{06}^{(2)}$ {\it south Hamiltonians}. Here, we will construct the transformation $r_3^{06}$.

It is known that the transformation $r_{06}$ can be obtained by doing ${\mathbb P}^2$-flip in the boundary of the variables $p_1,p_2$.

{\bf Step 1:} We make a change of variables:
$$
X^{(1)}=q_1 p_1^2, \quad Y^{(1)}=p_1, \quad Z^{(1)}=q_2 p_1^2, \quad W^{(1)}=\frac{p_2}{p_1}.
$$

{\bf Step 2:} We blow up along the curve $\{(X^{(1)},Y^{(1)},Z^{(1)},W^{(1)})|X^{(1)}=-Z^{(1)}+t, \ Y^{(1)}=0, \  W^{(1)}=-1\}$:
$$
X^{(2)}=\frac{X^{(1)}+Z^{(1)}-t}{Y^{(1)}}, \quad Y^{(2)}=Y^{(1)}, \quad Z^{(2)}=Z^{(1)}, \quad W^{(2)}=\frac{W^{(1)}+1}{Y^{(1)}}.
$$

{\bf Step 3:} We blow up along the surface $\{(X^{(2)},Y^{(2)},Z^{(2)},W^{(2)})|X^{(2)}=-2Z^{(2)} W^{(2)}+\alpha_0-\alpha_3+\alpha_6, \ Y^{(2)}=0 \}$:
$$
X^{(3)}=\frac{X^{(2)}+2Z^{(2)} W^{(2)}-\alpha_0+\alpha_3-\alpha_6}{Y^{(2)}}, \quad Y^{(3)}=Y^{(2)}, \quad Z^{(3)}=Z^{(2)}, \quad W^{(3)}=W^{(2)}.
$$
By taking the coordinate system as
$$
(x_3,y_3,z_3,w_3)=(X^{(3)},Y^{(3)},Z^{(3)},W^{(3)}),
$$
we can obtain the coordinate system $r_3^{06}$.

We note that the coordinate system $(X^{(1)},Y^{(1)},Z^{(1)},W^{(1)})$ can be obtained by the following transformations;

{\bf Step 1:} We take the following coordinate system (see Figure 4):
$$
x^{(1)}=\frac{1}{q_1}, \quad y^{(1)}=p_1, \quad z^{(1)}=\frac{q_2}{q_1}, \quad w^{(1)}=p_2.
$$

{\bf Step 2:} We blow up along the surface $\{(x^{(1)},y^{(1)},z^{(1)},w^{(1)})|x^{(1)}=y^{(1)}=0 \}$:
$$
x^{(2)}=\frac{x^{(1)}}{y^{(1)}}, \quad y^{(2)}=y^{(1)}, \quad z^{(2)}=z^{(1)}, \quad w^{(2)}=w^{(1)}.
$$

{\bf Step 3:} We blow up along the surface $\{(x^{(2)},y^{(2)},z^{(2)},w^{(2)})|x^{(2)}=y^{(2)}=0 \}$:
$$
x^{(3)}=\frac{x^{(2)}}{y^{(2)}}, \quad y^{(3)}=y^{(2)}, \quad z^{(3)}=z^{(2)}, \quad w^{(3)}=w^{(2)}.
$$

{\bf Step 4:} We blow up along the surface $\{(x^{(3)},y^{(3)},z^{(3)},w^{(3)})|x^{(3)}=z^{(3)}=0 \}$:
$$
x^{(4)}=x^{(3)}, \quad y^{(4)}=y^{(3)}, \quad z^{(4)}=\frac{z^{(3)}}{x^{(3)}}, \quad w^{(4)}=w^{(3)}.
$$

{\bf Step 5:} We make a change of variables (see \cite{MMT}):
$$
x^{(5)}=\frac{1}{x^{(4)}}, \quad y^{(5)}=y^{(4)}, \quad z^{(5)}=z^{(4)}, \quad w^{(5)}=w^{(4)}.
$$
We see that
$$
(X^{(1)},Y^{(1)},Z^{(1)},W^{(1)})=(x^{(5)},y^{(5)},z^{(5)},w^{(5)})=(q_1 p_1^2,p_1,q_2 p_1^2,p_2).
$$

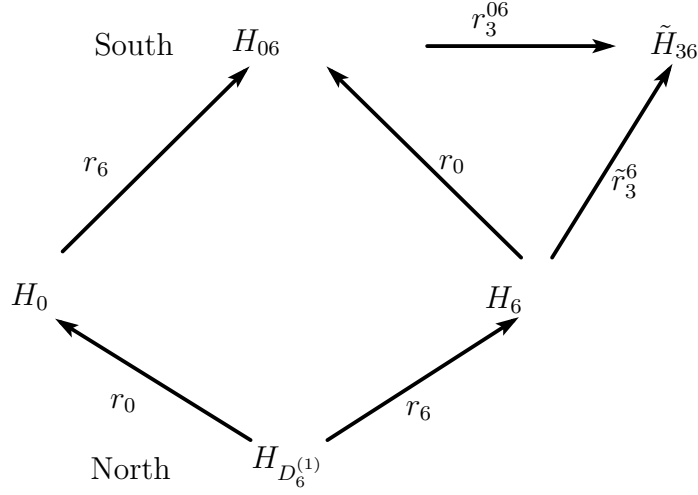
\begin{figure}
\unitlength 0.1in
\begin{picture}( 34.4600, 23.6500)( 31.6000,-31.4000)
\put(44.2000,-33.1000){\makebox(0,0)[lb]{$H_{D_6^{(1)}}$}}%
%
\special{pn 20}%
\special{pa 4416 3064}%
\special{pa 3422 2446}%
\special{fp}%
\special{sh 1}%
\special{pa 3422 2446}%
\special{pa 3468 2498}%
\special{pa 3466 2474}%
\special{pa 3488 2464}%
\special{pa 3422 2446}%
\special{fp}%
%
\special{pn 20}%
\special{pa 4818 3064}%
\special{pa 5792 2438}%
\special{fp}%
\special{sh 1}%
\special{pa 5792 2438}%
\special{pa 5726 2456}%
\special{pa 5748 2466}%
\special{pa 5748 2490}%
\special{pa 5792 2438}%
\special{fp}%
\put(31.6000,-23.8100){\makebox(0,0)[lb]{$H_0$}}%
\put(56.4200,-23.9800){\makebox(0,0)[lb]{$H_6$}}%
\put(43.1600,-10.4900){\makebox(0,0)[lb]{$H_{06}$}}%
%
\special{pn 20}%
\special{pa 3432 2076}%
\special{pa 4386 1130}%
\special{fp}%
\special{sh 1}%
\special{pa 4386 1130}%
\special{pa 4324 1162}%
\special{pa 4348 1168}%
\special{pa 4352 1190}%
\special{pa 4386 1130}%
\special{fp}%
%
\special{pn 20}%
\special{pa 5832 2108}%
\special{pa 4838 1130}%
\special{fp}%
\special{sh 1}%
\special{pa 4838 1130}%
\special{pa 4872 1190}%
\special{pa 4876 1166}%
\special{pa 4900 1162}%
\special{pa 4838 1130}%
\special{fp}%
\put(65.0500,-10.6500){\makebox(0,0)[lb]{$\tilde{H}_{36}$}}%
%
\special{pn 20}%
\special{pa 5994 2108}%
\special{pa 6606 1122}%
\special{fp}%
\special{sh 1}%
\special{pa 6606 1122}%
\special{pa 6554 1168}%
\special{pa 6578 1166}%
\special{pa 6588 1188}%
\special{pa 6606 1122}%
\special{fp}%
\put(36.8200,-29.1900){\makebox(0,0)[lb]{$r_0$}}%
\put(52.2900,-29.5900){\makebox(0,0)[lb]{$r_6$}}%
\put(35.4000,-16.8000){\makebox(0,0)[lb]{$r_6$}}%
\put(54.0100,-16.6700){\makebox(0,0)[lb]{$r_0$}}%
\put(55.7100,-9.4500){\makebox(0,0)[lb]{$r_3^{06}$}}%
\put(63.0400,-17.9500){\makebox(0,0)[lb]{$\tilde{r}_3^{6}$}}%
%
\special{pn 20}%
\special{pa 5340 1000}%
\special{pa 6294 1000}%
\special{fp}%
\special{sh 1}%
\special{pa 6294 1000}%
\special{pa 6228 980}%
\special{pa 6242 1000}%
\special{pa 6228 1020}%
\special{pa 6294 1000}%
\special{fp}%
\put(35.7800,-32.7100){\makebox(0,0)[lb]{North}}%
\put(35.9800,-10.4700){\makebox(0,0)[lb]{South}}%
\end{picture}%
\label{fig:•Ïg‹Èü'ÌŠG}
\caption{Relation between Hamiltonians $H_{D_6^{(1)}}$ and $H_0,H_6,H_{06},\tilde{H}_{36}$}
\end{figure}

Finally, we construct the holomorphy conditions of the Hamiltonian system with polynomial Hamiltonians $\tilde{H}_{36}^{(1)}=\tilde{r}_{3}^{6} (H_6^{(1)}), \  \tilde{H}_{36}^{(2)}=\tilde{r}_{3}^{6} (H_6^{(2)}),\\
\tilde{r}_3^{6}:x=-\left( \left(q_1+\frac{q_2}{p_1^2} \right)p_1-\frac{2q_2(p_2p_1+1)}{p_1}-(\alpha_3-\alpha_6) \right)p_1, \ y=\frac{1}{p_1},z=\frac{q_2}{p_1^2}, \ w=(p_2p_1+1)p_1$. These conditions are seven birational and symplectic transformations explicitly given by $R_i^{36} \ (i=0,1,\ldots,6)${\rm : \rm}

{\footnotesize
\begin{align*}
R_0^{36}:&x_0=q_1-\frac{2q_2p_2-\alpha_0+\alpha_3-\alpha_6}{p_1}+\frac{t}{p_1^2}, \ y_0=p_1, \ z_0=\frac{q_2}{p_1^2}, \ w_0=p_2p_1^2, \\
R_1^{36}:&x_1=q_1-\frac{2q_2p_2-\alpha_1+\alpha_3-\alpha_6}{p_1}+\frac{s}{p_1^2}, \ y_1=p_1, \ z_1=\frac{q_2}{p_1^2}, \ w_1=p_2p_1^2, \\
R_2^{36}:&x_2=-\left( \left(q_1-\frac{q_2}{p_1^2} \right)p_1-\frac{2q_2(p_2p_1-1)}{p_1}-(\alpha_2+\alpha_3-\alpha_6) \right)p_1, \ y_2=\frac{1}{p_1},\\
&z_2=\frac{q_2}{p_1^2} \ w_2=(p_2p_1-1)p_1, \\
R_3^{36}:&x_3=-\left( \left(q_1-\frac{q_2}{p_1^2} \right)p_1-\frac{2q_2(p_2p_1-1)}{p_1}-(\alpha_3-\alpha_6) \right)p_1, \ y_3=\frac{1}{p_1},\\
&z_3=\frac{q_2}{p_1^2}, \ w_3=(p_2p_1-1)p_1, \\
R_4^{36}:&x_4=q_1, \ y_4=p_1, \ z_4=-(q_2p_2-\alpha_4-\alpha_6)p_2, \ w_4=\frac{1}{p_2}, \\
R_5^{36}:&x_5=q_1+\frac{(\alpha_5-\alpha_6)p_2}{p_1p_2-1}+\frac{p_2^2}{(p_1p_2-1)^2}, \ y_5=p_1,\\
&z_5=q_2+\frac{(\alpha_5-\alpha_6)p_1}{p_1p_2-1}+\frac{1}{(p_1p_2-1)^2}, \ w_5=p_2, \\
R_6^{36}:&x_6=q_1, \ y_6=p_1, \ z_6=-(q_2p_2-\alpha_6)p_2, \ w_6=\frac{1}{p_2},
\end{align*}
where $R_0^{36} \left(\tilde{H}_{36}^{(1)}-\frac{1}{p_1} \right), \ R_1^{36} \left(\tilde{H}_{36}^{(2)}-\frac{1}{p_1} \right)$.
}

\vspace{0.5cm}

{\bf Hamiltonians $H_{03}^{(1)}=r_{03}(H_1^{D_6^{(1)}}-(p_1+p_2)), \ H_{03}^{(2)}=r_{03}(H_2^{D_6^{(1)}}), \ r_{03}:x=-((q_1-t)(p_1+p_2)-\alpha_0)(p_1+p_2), \ y=\frac{1}{p_1+p_2}, \ z=-((q_2-q_1)p_2-\alpha_3)p_2, \ w=\frac{1}{p_2}$}
\begin{align*}
&r_0^{03}:x_0=-(q_1p_1-\alpha_0)p_1, \ y_0=\frac{1}{p_1}, \ z_0=q_2, \ w_0=p_2, \\
&r_1^{03}:x_1=q_1+\frac{\alpha_1-\alpha_0}{p_1}+\frac{s-t}{p_1^2}, \ y_1=p_1, \ z_1=q_2, \ w_1=p_2, \\
&r_2^{03}:x_2=\frac{1}{q_1}, \ y_2=-((p_1-p_2)q_1+\alpha_2)q_1, \ z_2=q_2+q_1 \ w_2=p_2, \\
&r_3^{03}:x_3=q_1, \ y_3=p_1, \ z_3=-(q_2p_2-\alpha_3)p_2, \ w_3=\frac{1}{p_2}, \\
&r_4^{03}:x_4=q_1, \ y_4=p_1, \ z_4=-(q_2p_2-\alpha_4-\alpha_3)p_2, \ w_4=\frac{1}{p_2}, \\
&r_5^{03}:x_5=q_1 p_2^2, \ y_5=\frac{p_1-p_2}{p_2^2}, \ z_5=q_2-q_1+\frac{2q_1p_1+\alpha_5-\alpha_0-\alpha_3}{p_2}+\frac{1-t}{p_2^2}, \ w_5=p_2, \\
&r_6^{03}:x_6=q_1 p_2^2, \ y_6=\frac{p_1-p_2}{p_2^2}, \ z_6=q_2-q_1+\frac{2q_1p_1+\alpha_6-\alpha_0-\alpha_3}{p_2}-\frac{t}{p_2^2}, \ w_6=p_2,
\end{align*}
where $r_1^{03} \left(H_{03}^{(1)}+\frac{1}{p_1} \right), \ r_5^{03} \left(H_{03}^{(1)}+\frac{1}{p_2} \right), \ r_6^{03} \left(H_{03}^{(1)}+\frac{1}{p_2} \right), \ r_1^{03} \left(H_{03}^{(2)}-\frac{1}{p_1} \right)$.

\begin{figure}
\unitlength 0.1in
\begin{picture}( 49.6700, 32.6400)(  7.6000,-34.3100)
%
\special{pn 8}%
\special{pa 3418 368}%
\special{pa 1556 1396}%
\special{fp}%
\special{pa 1574 1380}%
\special{pa 4952 1380}%
\special{fp}%
\special{pa 3426 368}%
\special{pa 4932 1380}%
\special{fp}%
%
\special{pn 8}%
\special{pa 1594 1380}%
\special{pa 994 3192}%
\special{fp}%
\special{pa 994 3192}%
\special{pa 5644 3192}%
\special{fp}%
\special{pa 4932 1388}%
\special{pa 5644 3192}%
\special{fp}%
%
\special{pn 20}%
\special{pa 1014 3192}%
\special{pa 5624 3192}%
\special{fp}%
%
\special{pn 20}%
\special{pa 3436 368}%
\special{pa 3680 278}%
\special{fp}%
\special{sh 1}%
\special{pa 3680 278}%
\special{pa 3610 282}%
\special{pa 3630 296}%
\special{pa 3624 320}%
\special{pa 3680 278}%
\special{fp}%
%
\special{pn 20}%
\special{pa 1004 3192}%
\special{pa 920 3432}%
\special{fp}%
\special{sh 1}%
\special{pa 920 3432}%
\special{pa 960 3376}%
\special{pa 938 3382}%
\special{pa 922 3362}%
\special{pa 920 3432}%
\special{fp}%
%
\special{pn 20}%
\special{pa 5644 3192}%
\special{pa 5708 3432}%
\special{fp}%
\special{sh 1}%
\special{pa 5708 3432}%
\special{pa 5710 3362}%
\special{pa 5694 3380}%
\special{pa 5672 3372}%
\special{pa 5708 3432}%
\special{fp}%
%
\special{pn 20}%
\special{pa 3440 2240}%
\special{pa 3712 2136}%
\special{fp}%
\special{sh 1}%
\special{pa 3712 2136}%
\special{pa 3642 2142}%
\special{pa 3662 2156}%
\special{pa 3656 2180}%
\special{pa 3712 2136}%
\special{fp}%
\put(37.4400,-3.3700){\makebox(0,0)[lb]{$\frac{1}{p_1}$}}%
\put(37.5000,-22.9000){\makebox(0,0)[lb]{$\frac{1}{p_2}$}}%
\put(57.2700,-34.9800){\makebox(0,0)[lb]{$\frac{1}{q_2}$}}%
\put(7.6000,-35.5000){\makebox(0,0)[lb]{$\frac{1}{q_1}$}}%
\put(29.0000,-34.0000){\makebox(0,0)[lb]{$r_0^{36},r_1^{36}$}}%
%
\special{pn 8}%
\special{pa 1590 1400}%
\special{pa 2210 2010}%
\special{dt 0.045}%
\special{pa 2210 2010}%
\special{pa 2010 2600}%
\special{dt 0.045}%
\special{pa 2010 2600}%
\special{pa 990 3190}%
\special{dt 0.045}%
%
\special{pn 8}%
\special{pa 4930 1390}%
\special{pa 4400 2010}%
\special{dt 0.045}%
\special{pa 4400 2010}%
\special{pa 4610 2630}%
\special{dt 0.045}%
\special{pa 4610 2630}%
\special{pa 5620 3180}%
\special{dt 0.045}%
%
\special{pn 8}%
\special{pa 2200 2020}%
\special{pa 4380 2020}%
\special{dt 0.045}%
%
\special{pn 8}%
\special{pa 2030 2610}%
\special{pa 4620 2630}%
\special{dt 0.045}%
%
\special{pn 8}%
\special{pa 2220 2000}%
\special{pa 3420 1470}%
\special{dt 0.045}%
\special{pa 3420 1470}%
\special{pa 3420 1470}%
\special{dt 0.045}%
%
\special{pn 8}%
\special{pa 3430 1480}%
\special{pa 4390 2010}%
\special{dt 0.045}%
%
\special{pn 8}%
\special{pa 3450 390}%
\special{pa 3450 2230}%
\special{dt 0.045}%
%
\special{pn 8}%
\special{pa 3440 2230}%
\special{pa 2020 2610}%
\special{dt 0.045}%
%
\special{pn 8}%
\special{pa 3450 2240}%
\special{pa 4600 2620}%
\special{dt 0.045}%
%
\special{pn 20}%
\special{pa 1590 1390}%
\special{pa 2210 2010}%
\special{fp}%
%
\special{pn 20}%
\special{pa 3070 2020}%
\special{pa 3570 3170}%
\special{fp}%
%
\special{pn 20}%
\special{pa 4930 1400}%
\special{pa 5640 3180}%
\special{fp}%
\put(18.7000,-17.2000){\makebox(0,0)[lb]{$r_2^{36}$}}%
\put(34.8000,-29.5000){\makebox(0,0)[lb]{$r_4^{36}$}}%
\put(54.1000,-25.2000){\makebox(0,0)[lb]{$r_5^{36}$}}%
\put(32.2000,-10.5000){\makebox(0,0)[lb]{$r_3^{36}$}}%
\put(32.6000,-20.0000){\makebox(0,0)[lb]{$r_6^{36}$}}%
\end{picture}%
\label{fig:D6r36'Ì"ÁˆÙ"_2}
\caption{Accessible singular loci of the system with Hamiltonians $H_{36}^{(1)},H_{36}^{(2)}$}
\end{figure}
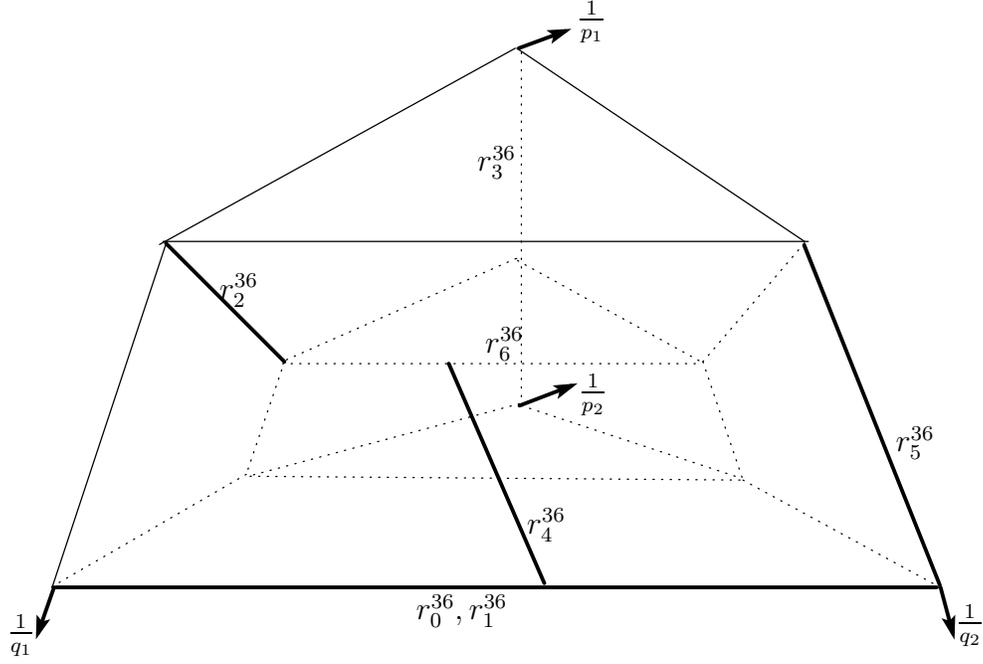

{\bf Hamiltonians $H_{36}^{(1)}=r_{36}(H_1^{D_6^{(1)}}), \ H_{36}^{(2)}=r_{36}(H_2^{D_6^{(1)}}), \ r_{36}:x=-((q_1-q_2)p_1-\alpha_3)p_1, \ y=\frac{1}{p_1}, \ z=-(q_2(p_2+p_1)-\alpha_6)(p_2+p_1), \ w=\frac{1}{p_2+p_1}$}
\begin{align*}
&r_0^{36}:x_0=q_1p_2^2, \ y_0=\frac{p_1-p_2}{p_2^2}, \ z_0=q_2-q_1+\frac{2q_1p_1+\alpha_0-\alpha_3-\alpha_6}{p_2}+\frac{t}{p_2^2}, \ w_0=p_2, \\
&r_1^{36}:x_1=q_1p_2^2, \ y_1=\frac{p_1-p_2}{p_2^2}, \ z_1=q_2-q_1+\frac{2q_1p_1+\alpha_1-\alpha_3-\alpha_6}{p_2}+\frac{s}{p_2^2}, \ w_1=p_2, \\
&r_2^{36}:x_2=-(q_1p_1-\alpha_2-\alpha_3)p_1, \ y_2=\frac{1}{p_1}, \ z_2=q_2 \ w_2=p_2, \\
&r_3^{36}:x_3=-(q_1p_1-\alpha_3)p_1, \ y_3=\frac{1}{p_1}, \ z_3=q_2, \ w_3=p_2, \\
&r_4^{36}:x_4=\frac{1}{q_1}, \ y_4=-((p_1-p_2)q_1+\alpha_4)q_1, \ z_4=q_2+q_1, \ w_4=p_2, \\
&r_5^{36}:x_5=q_1, \ y_5=p_1, \ z_5=q_2+\frac{\alpha_5-\alpha_6}{p_2}+\frac{1}{p_2^2}, \ w_5=p_2, \\
&r_6^{36}:x_6=q_1, \ y_6=p_1, \ z_6=-(q_2p_2-\alpha_6)p_2, \ w_6=\frac{1}{p_2},
\end{align*}
where $r_0^{36} \left(H_{36}^{(1)}-\frac{1}{p_2} \right), \ r_1^{36} \left(H_{36}^{(2)}-\frac{1}{p_2} \right)$.

By doing augmentation theory by N. Tahara (see \cite{Ta}), we can also obtain the Hamiltonians $r_4^{36}(H_{36}^{(1)}),r_4^{36}(H_{36}^{(2)})$.

\vspace{0.5cm}

\begin{figure}
\unitlength 0.1in
\begin{picture}( 47.3000, 42.9000)( 19.3000,-56.1000)
\put(43.8000,-40.2500){\makebox(0,0)[lb]{$H_{D_6^{(1)}}$}}%
%
\special{pn 20}%
\special{pa 4376 3778}%
\special{pa 3382 3160}%
\special{fp}%
\special{sh 1}%
\special{pa 3382 3160}%
\special{pa 3428 3212}%
\special{pa 3426 3188}%
\special{pa 3448 3178}%
\special{pa 3382 3160}%
\special{fp}%
%
\special{pn 20}%
\special{pa 4778 3778}%
\special{pa 5752 3152}%
\special{fp}%
\special{sh 1}%
\special{pa 5752 3152}%
\special{pa 5686 3172}%
\special{pa 5708 3182}%
\special{pa 5708 3206}%
\special{pa 5752 3152}%
\special{fp}%
\put(31.2000,-30.9600){\makebox(0,0)[lb]{$H_3$}}%
\put(56.0200,-31.1300){\makebox(0,0)[lb]{$H_6$}}%
\put(44.9000,-17.7000){\makebox(0,0)[lb]{$H_{36}$}}%
%
\special{pn 20}%
\special{pa 3392 2792}%
\special{pa 4346 1844}%
\special{fp}%
\special{sh 1}%
\special{pa 4346 1844}%
\special{pa 4284 1878}%
\special{pa 4308 1882}%
\special{pa 4312 1906}%
\special{pa 4346 1844}%
\special{fp}%
\put(64.6500,-17.8000){\makebox(0,0)[lb]{$\tilde{H}_{36}$}}%
\put(36.4200,-36.3400){\makebox(0,0)[lb]{$r_3$}}%
\put(51.8900,-36.7400){\makebox(0,0)[lb]{$r_6$}}%
\put(35.0000,-23.9500){\makebox(0,0)[lb]{$r_6$}}%
\put(55.3100,-16.6000){\makebox(0,0)[lb]{$r_0^{36}$}}%
\put(66.4000,-24.0000){\makebox(0,0)[lb]{$R_6^{0}$}}%
%
\special{pn 20}%
\special{pa 5300 1716}%
\special{pa 6254 1716}%
\special{fp}%
\special{sh 1}%
\special{pa 6254 1716}%
\special{pa 6188 1696}%
\special{pa 6202 1716}%
\special{pa 6188 1736}%
\special{pa 6254 1716}%
\special{fp}%
\put(43.8000,-43.0000){\makebox(0,0)[lb]{North}}%
\put(43.7000,-14.9000){\makebox(0,0)[lb]{South}}%
%
\special{pn 20}%
\special{pa 4130 1720}%
\special{pa 3380 1720}%
\special{fp}%
\special{sh 1}%
\special{pa 3380 1720}%
\special{pa 3448 1740}%
\special{pa 3434 1720}%
\special{pa 3448 1700}%
\special{pa 3380 1720}%
\special{fp}%
\put(19.3000,-17.7000){\makebox(0,0)[lb]{$H_{364}=r_4^{36}(H_{36})$}}%
\put(31.2000,-20.2000){\makebox(0,0)[lb]{Augmentation}}%
%
\special{pn 20}%
\special{pa 6590 1880}%
\special{pa 6590 2900}%
\special{fp}%
\special{sh 1}%
\special{pa 6590 2900}%
\special{pa 6610 2834}%
\special{pa 6590 2848}%
\special{pa 6570 2834}%
\special{pa 6590 2900}%
\special{fp}%
\put(64.0000,-31.2000){\makebox(0,0)[lb]{$H'_6$}}%
\put(65.0000,-40.3000){\makebox(0,0)[lb]{$\tilde{H}_{D_6^{(1)}}$}}%
\put(66.6000,-35.7000){\makebox(0,0)[lb]{$r_6$}}%
\put(58.9000,-38.7000){\makebox(0,0)[lb]{$\varphi$}}%
%
\special{pn 20}%
\special{pa 6590 3160}%
\special{pa 6590 3650}%
\special{fp}%
\special{sh 1}%
\special{pa 6590 3650}%
\special{pa 6610 3584}%
\special{pa 6590 3598}%
\special{pa 6570 3584}%
\special{pa 6590 3650}%
\special{fp}%
%
\special{pn 20}%
\special{pa 6400 3970}%
\special{pa 5520 3970}%
\special{fp}%
\special{sh 1}%
\special{pa 5520 3970}%
\special{pa 5588 3990}%
\special{pa 5574 3970}%
\special{pa 5588 3950}%
\special{pa 5520 3970}%
\special{fp}%
\put(24.3000,-51.8000){\makebox(0,0)[lb]{$R_6^{0}:x_3=-\left( \left(q_1+\frac{q_2}{p_1^2} \right)p_1-\frac{2q_2(p_2p_1+1)}{p_1}-(\alpha_0-\alpha_6) \right)p_1, \ y_3=\frac{1}{p_1},$}}%
\put(24.4000,-57.8000){\makebox(0,0)[lb]{$\varphi:x=q_1+q_2+t, \ y=p_1, \ z=q_2, \ w=p_2-p_1$}}%
\put(31.6000,-54.6000){\makebox(0,0)[lb]{$z_3=\frac{q_2}{p_1^2}, \ w_3=(p_2p_1+1)p_1$}}%
\put(24.3000,-45.8000){\makebox(0,0)[lb]{$r_0^{36}:x_0=q_1-q_2+\frac{2q_2 p_2+\alpha_0-\alpha_3-\alpha_6}{p_1}+\frac{t}{p_1^2}, \ y_0=p_1,$}}%
\put(31.6000,-48.6000){\makebox(0,0)[lb]{$z_0=q_2 p_1^2, \ w_0=\frac{p_2-p_1}{p_1^2}$}}%
\end{picture}%
\label{fig:•Ïg‹Èü'Æ"cŒ´—˜_'ÌŠG}
\caption{Relation between Hamiltonians $H_{D_6^{(1)}}$ and $H_3,H_6,H_{36},\tilde{H}_{36},H_{364}$}
\end{figure}
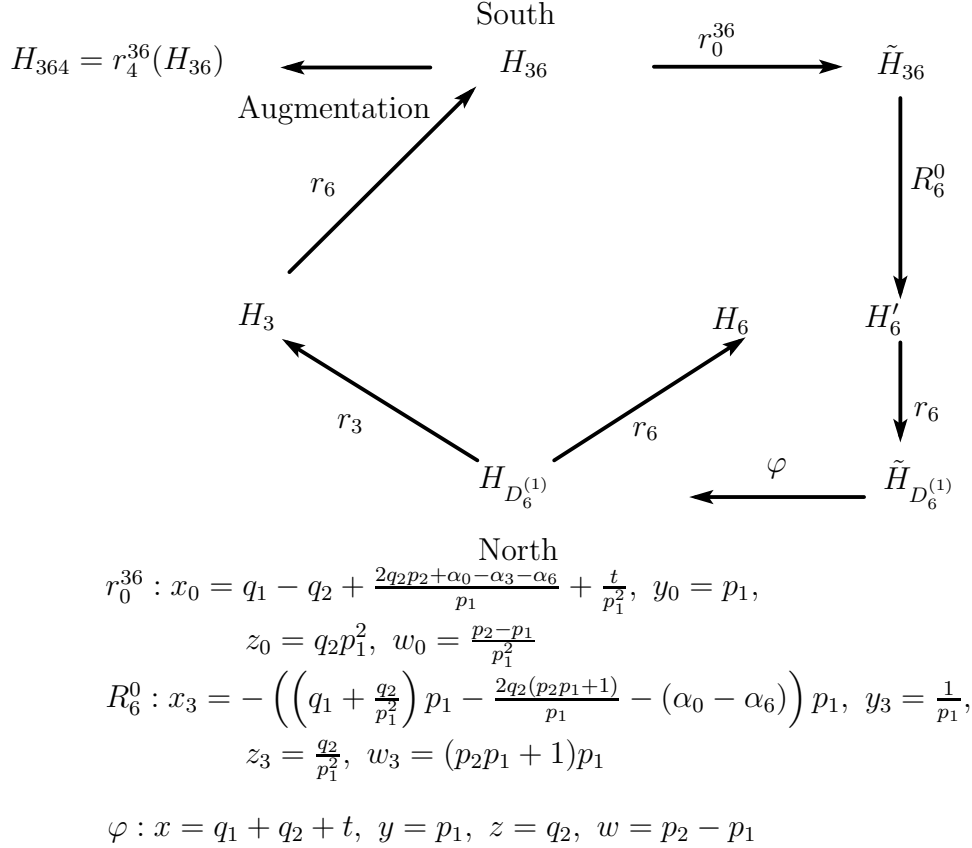

{\bf Hamiltonians $H_{364}^{(1)}=r_4^{36}(H_{36}^{(1)}), \ H_{364}^{(2)}=r_4^{36}(H_{36}^{(2)}), \ r_4^{36}:x_4=\frac{1}{q_1}, \ y_4=-((p_1-p_2)q_1+\alpha_4)q_1, \ z_4=q_2+q_1, \ w_4=p_2$}
\begin{align*}
&r_0^{364}:x_0=x_4 w_4^2, \ y_0=\frac{y_4}{w_4^2}, \ z_0=z_4+\frac{2x_4y_4+\alpha_0-(\alpha_3+\alpha_6)}{w_4}+\frac{t}{w_4^2}, \ w_0=w_4, \\
&r_1^{364}:x_1=x_4 w_4^2, \ y_1=\frac{y_4}{w_4^2}, \ z_1=z_4+\frac{2x_4y_4+\alpha_1-(\alpha_3+\alpha_6)}{w_4}+\frac{s}{w_4^2}, \ w_1=w_4, \\
&r_2^{364}:x_2=\frac{1}{q_1}, \ y_2=-(p_1q_1+\alpha_2+\alpha_3+\alpha_4)q_1, \ z_2=q_2, \ w_2=p_2, \\
&r_3^{364}:x_3=\frac{1}{q_1}, \ y_3=-(p_1q_1+\alpha_3+\alpha_4)q_1, \ z_3=q_2, \ w_3=p_2, \\
&r_4^{364}:x_4=\frac{1}{q_1}, \ y_4=-(p_1q_1+\alpha_4)q_1, \ z_4=q_2, \ w_4=p_2, \\
&r_5^{364}:x_5=q_1, \ y_5=p_1, \ z_5=q_2+\frac{\alpha_5-\alpha_6}{p_2}+\frac{1}{p_2^2}, \ w_5=p_2, \\
&r_6^{364}:x_6=q_1, \ y_6=p_1-\frac{\alpha_6 q_2}{q_1 q_2-1}, \ z_6=q_2, \ w_6=p_2-\frac{\alpha_6 q_1}{q_1 q_2-1},
\end{align*}
where $r_0^{364} \left(r_4^{364}(H_{364}^{(1)})-\frac{1}{w_4} \right), \ r_1^{364} \left(r_4^{364}(H_{364}^{(2)})-\frac{1}{w_4} \right)$.

For the Hamiltonians $H_{35}^{(1)}=r_{35}(H_1^{D_6^{(1)}}), \ H_{35}^{(2)}=r_{35}(H_2^{D_6^{(1)}}), \ r_{35}:x=-((q_1-q_2)p_1-\alpha_3)p_1, \ y=\frac{1}{p_1}, \ z=-((q_2-1)(p_2+p_1)-\alpha_5)(p_2+p_1), \ w=\frac{1}{p_2+p_1}$, we can obtain its holomorphy conditions by the same way.

\begin{center}
\begin{tabular}{|c||c|} \hline
 & Holomorphy conditions \\ \hline
$r_3$ & $x_3=-((q_1-q_2)p_1-\alpha_3)p_1, \ y_3=\frac{1}{p_1}, \ z_3=q_2, \ w_3=p_2+p_1$  \\ \hline
$\tilde{r}_{3}$ & $\tilde{x}_3=q_1, \ \tilde{y}_3=p_1+p_2, \ \tilde{z}_3=-((q_2-q_1)p_2-\alpha_3)p_2, \ \tilde{w}_3=\frac{1}{p_2}$  \\ \hline
$r_3^{0}$ & $x_3=q_1+\frac{\alpha_3-\alpha_0}{p_1}+\frac{q_2-t}{p_1^2}, \ y_3=p_1, \ z_3=q_2, \ w_3=p_2+\frac{1}{p_1}$  \\ \hline
$\tilde{r}_3^{6}$ & $x_3=-\left( \left(q_1+\frac{q_2}{p_1^2} \right)p_1-\frac{2q_2(p_2p_1+1)}{p_1}-(\alpha_3-\alpha_6) \right)p_1, \ y_3=\frac{1}{p_1},z_3=\frac{q_2}{p_1^2}, \ w_3=(p_2p_1+1)p_1$  \\ \hline
$r_3^{06}$ & $x_3=q_1+q_2+\frac{2q_2p_2+\alpha_3-\alpha_0-\alpha_6}{p_1}-\frac{t}{p_1^2}, \ y_3=p_1, \ z_3=q_2p_1^2, \ w_3=\frac{p_2+p_1}{p_1^2}$  \\ \hline
\end{tabular}
\end{center}
In the above table, we summarize some holomorphy conditions associated with the transformation $r_3$.

\section{Appendix D: Fuji-Suzuki system}

{\bf Holomorphy conditions}

Define birational and symplectic transformations $r_i \ (i=0,1,\ldots,6)$ as follows:

\begin{align}\label{holo;FS}
\begin{split}
r_0:(x_0,y_0,z_0,w_0)=&\left(-((q_1-q_2)p_1-\beta_0)p_1,\frac{1}{p_1},q_2,p_2+p_1 \right),\\
r_1:(x_1,y_1,z_1,w_1)=&\left(\frac{1}{q_1},-(q_1p_1+\beta_1)q_1,q_2,p_2 \right),\\
r_2:(x_2,y_2,z_2,w_2)=&\left(-((q_1-t)p_1-\beta_2)p_1,\frac{1}{p_1},q_2,p_2 \right),\\
r_3:(x_3,y_3,z_3,w_3)=&\left(-(q_1p_1+q_2p_2-\beta_3)p_1,\frac{1}{p_1},q_2p_1,\frac{p_2}{p_1} \right),\\
r_4:(x_4,y_4,z_4,w_4)=&\left(q_1,p_1,-((q_2-s)p_2-\beta_4)p_2,\frac{1}{p_2} \right),\\
r_5:(x_5,y_5,z_5,w_5)=&\left(q_1,p_1,\frac{1}{q_2},-(p_2q_2+\beta_5)q_2 \right),\\
r_6:(x_6,y_6,z_6,w_6)=&\left(-((q_1-1)p_1+(q_2-1)p_2-\beta_6)p_1,\frac{1}{p_1},(q_2-1)p_1,\frac{p_2}{p_1} \right).
\end{split}
\end{align}
There exist two polynomials $H_1^{FS}$ and $H_2^{FS}$, such that the Hamiltonian system
\begin{equation}\label{eq:FS}
  \left\{
  \begin{aligned}
   dq_1 =&\frac{\partial H_1^{FS}}{\partial p_1}dt+\frac{\partial H_2^{FS}}{\partial p_1}ds, \quad dp_1 =-\frac{\partial H_1^{FS}}{\partial q_1}dt-\frac{\partial H_2^{FS}}{\partial q_1}ds,\\
   dq_2 =&\frac{\partial H_1^{FS}}{\partial p_2}dt+\frac{\partial H_2^{FS}}{\partial p_2}ds, \quad dp_2 =-\frac{\partial H_1^{FS}}{\partial q_2}dt-\frac{\partial H_2^{FS}}{\partial q_2}ds
   \end{aligned}
  \right.
\end{equation}
is transformed into a polynomial Hamiltonian system under the action of each $r_i \ (i=0,1,\ldots,6)$, where two polynomial Hamiltonians $H_1^{FS},H_2^{FS}$ are given by {\rm (cf. \cite{FS1,KNS,Fuji}) \rm}
\begin{align}
\begin{split}
H_1^{FS} &=H_{VI}(q_1,p_1,t,s;\beta_2,\beta_0+\beta_4,\beta_1,\beta_5+\beta_6,\beta_3+\beta_5)\\
&+H_{VI}(q_2,p_2,t,s;\beta_0+\beta_2,\beta_4,\beta_5,\beta_1+\beta_6,\beta_1+\beta_3)\\
&+\frac{(q_1-t)(q_2-s) \{2(q_1 p_1+\beta_1)(q_2p_2+\beta_5)-(q_1 p_1+\beta_1)p_2-(q_2p_2+\beta_5)p_1 \}}{(\beta_0+2\beta_1+\beta_2+\beta_3+\beta_4+2\beta_5+\beta_6)t(t-1)(t-s)}\\
&-\frac{ 2\beta_1 \beta_5 s }{(\beta_0+2\beta_1+\beta_2+\beta_3+\beta_4+2\beta_5+\beta_6)(t-1)(t-s)},\\
H_2^{FS} &=\pi(H_1^{FS}), \quad \pi=\{q_1 \leftrightarrow q_2, \ p_1 \leftrightarrow p_2, \  t \leftrightarrow s, \ \beta_1 \leftrightarrow \beta_5, \ \beta_2 \leftrightarrow \beta_4 \}.
\end{split}
\end{align}
The symbol $H_{VI}(q,p,t,\eta;\gamma_0,\gamma_1,\gamma_2,\gamma_3,\gamma_4)$ denotes (see \cite{sasa3})
\begin{align}
\begin{split}
&t(t-1)(t-\eta) H_{VI}(q,p,t,\eta;\gamma_0,\gamma_1,\gamma_2,\gamma_3,\gamma_4)\\
&=q(q-1)(q-\eta)(q-t)p^2+\{\gamma_1(t-\eta)q(q-1)+2\gamma_2 q(q-1)(q-\eta)\\
&+\gamma_3 (t-1)q(q-\eta)+\gamma_4 t(q-1)(q-\eta) \}p\\
&+\gamma_2 \{(\gamma_1+\gamma_2)(t-\eta)+\gamma_2(q-1)+\gamma_3(t-1)+t \gamma_4 \}q, \quad (\gamma_0+\gamma_1+2\gamma_2+\gamma_3+\gamma_4=1).
\end{split}
\end{align}

\vspace{0.5cm}

We note that the holomorphy conditions should be read that in the Hamiltonian $H_1^{FS}$
\begin{align*}
\begin{split}
&r_2(H_1^{FS} - p_1)
\end{split}
\end{align*}
are polynomials with respect to $x_2,y_2,z_2,w_2$, and in the Hamiltonian $H_2^{FS}$
\begin{align*}
\begin{split}
&r_4(H_2^{FS} -p_2)
\end{split}
\end{align*}
are polynomials with respect to $x_4,y_4,z_4,w_4$.

We see that the birational and symplectic transformation $\varphi$\rm{:\rm}
\begin{equation}
  \left\{
  \begin{aligned}
   Q_1=&\frac{1-q_1}{q_1}, \quad P_1=-(p_1 q_1+\beta_1)q_1, \quad Q_2=\frac{1-q_2}{q_2}, \quad P_2=-(p_2 q_2+\beta_5)q_2,\\
   T=&\frac{1-t}{t}, \quad S=\frac{1-s}{s}
   \end{aligned}
  \right. 
\end{equation}
takes the system \eqref{eq:FS} into Fuji-Suzuki system (see \cite{FS1}) when $S=1$.

We remark that the relations between $\alpha_i \ (i=0,1,\ldots,5), \eta$ (see \cite{FS1}) and $\beta_j \ (j=0,1,\ldots,6)$ are explicitly given as follows:
\begin{align}\label{FSprelation}
\begin{split}
&\alpha_3=\beta_1+\beta_3+\beta_5+\beta_6, \quad \eta=\beta_1+\beta_3+\beta_5,\\
&\alpha_0=\beta_0, \ \alpha_1=\beta_1, \ \alpha_2=\beta_2, \ \alpha_4=\beta_4 , \alpha_5=\beta_5.
\end{split}
\end{align}
Of course, $\alpha_i$ and $\beta_j$ satisfy the relation:
\begin{align}
\begin{split}
&\alpha_0+\alpha_1+\alpha_2+\alpha_3+\alpha_4+\alpha_5=\beta_0+2\beta_1+\beta_2+\beta_3+\beta_4+2\beta_5+\beta_6=1.
\end{split}
\end{align}

\vspace{0.5cm}

{\bf Completely integrable}

{\footnotesize
\begin{proposition}
Setting
\begin{equation}
K_1:=-H_1+\frac{\beta_3(\beta_1+\beta_5)( \rm{Log \rm}(s-t)-\rm{Log \rm}(s-1) )}{(\beta_0+2\beta_1+\beta_2+\beta_3+\beta_4+2\beta_5+\beta_6)(t-1)^2}, \quad K_2:=-H_2.
\end{equation}
Two Hamiltonians $K_1$ and $K_2$ satisfy
\begin{equation}
\{K_1,K_2\}+\left(\frac{\partial}{\partial s} \right)K_1-\left(\frac{\partial}{\partial t} \right)K_2=0,
\end{equation}
where  $\{,\}$ denotes the Poisson brackets:
\begin{equation}
\{L_1,L_2\}=\frac{\partial L_1}{\partial p_1}\frac{\partial L_2}{\partial q_1}-\frac{\partial L_1}{\partial q_1}\frac{\partial L_2}{\partial p_1}+\frac{\partial L_1}{\partial p_2}\frac{\partial L_2}{\partial q_2}-\frac{\partial L_1}{\partial q_2}\frac{\partial L_2}{\partial p_2}.
\end{equation}
\end{proposition}
}

We remark that on new constant complex parameters $\beta_j \ (j=0,1,\ldots,6)$ the Hamiltonian system \eqref{eq:FS} is invariant under these birational and symplectic transformations $s_0,s_1,\ldots,s_9$ (cf. Appendix B in  \cite{FS1}), whose generators are defined as follows$:$ with {\it the notation} $(*):=(q_1,p_1,q_2,p_2,t,s;\beta_0,\beta_1,\ldots,\beta_6)$;

{\footnotesize

\begin{align}\label{symmetry:FS}
\begin{split}
&s_0:(*) \rightarrow \left(q_1,p_1-\frac{\beta_0}{q_1-q_2},q_2,p_2+\frac{\beta_0}{q_1-q_2},t,s;-\beta_0,\beta_1+\beta_0,\beta_2,\beta_3-\beta_0,\beta_4,\beta_5+\beta_0,\beta_6-\beta_0 \right),\\
&s_1:(*) \rightarrow \left(q_1+\frac{\beta_1}{p_1},p_1,q_2,p_2,t,s;\beta_0+\beta_1,-\beta_1,\beta_2+\beta_1,\beta_3+\beta_1,\beta_4,\beta_5,\beta_6+\beta_1 \right),\\
&s_2:(*) \rightarrow \left(q_1,p_1-\frac{\beta_2}{q_1-t},q_2,p_2,t,s;\beta_0,\beta_1+\beta_2,-\beta_2,\beta_3,\beta_4,\beta_5,\beta_6 \right),\\
&s_3:(*) \rightarrow (\frac{q_1 (-q_1 p_1 + q_1^2 p_1 - p_2 q_2 + p_2 q_2^2 - \beta_1 + 
    q_1 \beta_1 - \beta_5 + q_2 \beta_5 - \beta_6)}{g_1},\\
& -\frac{g_1 (q_1 p_1^2 -
        q_1^2 p_1^2 + p_2 p_1 q_2 - p_2 p_1 q_2^2 + \beta_1^2 - p_1 \beta_3 + 
       q_1 p_1 \beta_3 + \beta_1 \beta_3 + q_1 p_1 \beta_5 - 
       p_1 q_2 \beta_5 + \beta_1 \beta_5 + 
       q_1 p_1 \beta_6 + \beta_1 \beta_6)}{(-q_1 p_1 + q_1^2 p_1 - p_2 q_2 + 
      p_2 q_2^2 + q_1 \beta_1 + \beta_3 + q_2 \beta_5) (-q_1 p_1 + q_1^2 p_1 - 
      p_2 q_2 + p_2 q_2^2 - \beta_1 + q_1 \beta_1 - \beta_5 + 
      q_2 \beta_5 - \beta_6)},\\
& \frac{
 q_2 (-q_1 p_1 + q_1^2 p_1 - p_2 q_2 + p_2 q_2^2 - \beta_1 + q_1 \beta_1 - \beta_5 +
     q_2 \beta_5 - \beta_6)}{g_2},\\
& -\frac{g_2 (p_2 q_1 p_1 - 
       p_2 q_1^2 p_1 + p_2^2 q_2 - p_2^2 q_2^2 - p_2 q_1 \beta_1 + p_2 q_2 \beta_1 - 
       p_2 \beta_3 + 
       p_2 q_2 \beta_3 + \beta_1 \beta_5 + \beta_3 \beta_5 + \
\beta_5^2 + p_2 q_2 \beta_6 + \beta_5 \beta_6)}{(-q_1 p_1 + q_1^2 p_1 - 
      p_2 q_2 + p_2 q_2^2 + q_1 \beta_1 + \beta_3 + q_2 \beta_5) (-q_1 p_1 + 
      q_1^2 p_1 - p_2 q_2 + p_2 q_2^2 - \beta_1 + q_1 \beta_1 - \beta_5 + 
      q_2 \beta_5 - \beta_6)},\\
&t,s;\beta_0,\beta_1,\beta_1+\beta_2+\beta_3+\beta_5+\beta_6,-\beta_1-\beta_5-\beta_6,\beta_1+\beta_3+\beta_4+\beta_5+\beta_6,\beta_5,-\beta_1-\beta_3-\beta_5),\\
&s_4:(*) \rightarrow \left(q_1,p_1,q_2,p_2-\frac{\beta_4}{q_2-s},t,s;\beta_0,\beta_1,\beta_2,\beta_3,-\beta_4,\beta_5+\beta_4,\beta_6 \right),\\
&s_5:(*) \rightarrow \left(q_1,p_1,q_2+\frac{\beta_5}{p_2},p_2,t,s;\beta_0+\beta_5,\beta_1,\beta_2,\beta_3+\beta_5,\beta_4+\beta_5,-\beta_5,\beta_6+\beta_5 \right),\\
&s_6:(*) \rightarrow \left(1-q_1,-p_1,1-q_2,-p_2,1-t,1-s;\beta_0,\beta_1,\beta_2,\beta_6,\beta_4,\beta_5,\beta_3 \right),\\
&s_7:(*) \rightarrow \left(q_2,p_2,q_1,p_1,s,t;\beta_0,\beta_5,\beta_4,\beta_3,\beta_2,\beta_1,\beta_6 \right),\\
&s_8:(*) \rightarrow (1 - q_1, \frac{q_1 p_1 - q_1^2 p_1 + p_2 q_2 - p_2 q_2^2 - \beta_3 + q_1 \beta_3 + 
 q_1 \beta_5 - q_2 \beta_5 + 
 q_1 \beta_6}{(-1 + q_1) q_1}, \frac{(-1 + s) (-1 + q_1) q_2}{
s q_1 + q_2 - s q_2 - q_1 q_2}, \\
&-\frac{(s q_1 + q_2 - s q_2 - q_1 q_2) (-s p_2 q_1 - p_2 q_2 + s p_2 q_2 + 
    p_2 q_1 q_2 - \beta_5 + s \beta_5 + q_1 \beta_5)}{(-1 + s) s (-1 + 
    q_1) q_1},\\
&1-t,1-s;\beta_4,\beta_1+\beta_3+\beta_5+\beta_6,\beta_2,-\beta_5-\beta_6,\beta_0,\beta_5,-\beta_3-\beta_5 ),\\
&s_9:=s_7 \circ s_8 \quad ((s_{9})^6=1),
\end{split}
\end{align}
where 
\begin{align*}
\begin{split}
&g_1:=-q_1 p_1 + 
  q_1^2 p_1 - p_2 q_2 + p_2 q_2^2 + \beta_3 - q_1 \beta_3 - q_1 \beta_5 + 
  q_2 \beta_5 - 
  q_1 \beta_6,\\
&g_2:=-q_1 p_1 + q_1^2 p_1 - p_2 q_2 + p_2 q_2^2 + 
  q_1 \beta_1 - q_2 \beta_1 + \beta_3 - q_2 \beta_3 - 
  q_2 \beta_6.
\end{split}
\end{align*}
We note that the subgroup $<s_0,s_1,\ldots,s_5>$ generated by $s_0,s_1,\ldots,s_5$ is isomorphic to the affine Weyl group of type $A_5^{(1)}$  (see Appendix B in  \cite{FS1}), and the transformation $s_6$ was found by Professor K. Fuji in Kobe university in August 2012.

Finally, let us define the following translation operators{\rm : \rm}
\begin{align}
\begin{split}
&T_1:=(s_2s_9  s_9  s_1)^4, \quad T_2:=s_1T_1s_1,  \quad T_3:=s_5T_1s_5.
\end{split}
\end{align}
These translation operators act on parameters $\beta_i$ as follows$:$
\begin{align}
\begin{split}
T_1(\beta_0,\beta_1,\ldots,\beta_6)=&(\beta_0,\beta_1,\ldots,\beta_6)+(0,-1,1,0,-1,1,0),\\
T_2(\beta_0,\beta_1,\ldots,\beta_6)=&(\beta_0,\beta_1,\ldots,\beta_6)+(-1,1,0,-1,-1,1,-1),\\
T_3(\beta_0,\beta_1,\ldots,\beta_6)=&(\beta_0,\beta_1,\ldots,\beta_6)+(1,-1,1,1,0,-1,1).
\end{split}
\end{align}}

Finally, we remark some holomorphy conditions of the system \eqref{eq:FS}.

{\bf Hamiltonians $H_{04}^{(1)}=r_{04}(H_1^{FS}), \ H_{04}^{(2)}=r_{04}(H_2^{FS}-p_1-p_2), \ r_{04}:x=-((q_1-q_2)p_1-\beta_0)p_1, \ y=\frac{1}{p_1}, \ z=-((q_2-s)(p_2+p_1)-\beta_4)(p_2+p_1), \ w=\frac{1}{p_2+p_1}$}
\begin{align*}
&r_0^{04}:x_0=-(q_1p_1-\beta_0)p_1, \ y_0=\frac{1}{p_1}, \ z_0=q_2, \ w_0=p_2, \\
&r_1^{04}:x_1=-(q_1p_1-\beta_1-\beta_0)p_1, \ y_1=\frac{1}{p_1}, \ z_1=q_2, \ w_1=p_2, \\
&r_2^{04}:x_2=q_1-q_2+\frac{2q_1p_1+\beta_2-(\beta_0+\beta_4)}{p_1}+\frac{t-s}{p_1^2}, \ y_2=p_1, \ z_2=q_2 p_1^2 \ w_2=\frac{p_2-p_1}{p_1^2}, \\
&r_3^{04}:x_3=q_1p_2, \ y_3=\frac{p_1}{p_2}, \ z_3=q_2+\frac{q_1p_1+\beta_3-(\beta_0+\beta_4)}{p_2}-\frac{s}{p_2^2}, \ w_3=p_2, \\
&r_4^{04}:x_4=q_1, \ y_4=p_1, \ z_4=-(q_2 p_2-\beta_4)p_2, \ w_4=\frac{1}{p_2}, \\
&r_5^{04}:x_5=\frac{1}{q_1}, \ y_5=-((p_1-p_2)q_1+\beta_5)q_1, \ z_5=q_2+q_1, \ w_5=p_2,\\
&r_6^{04}:x_6=q_1p_2, \ y_6=\frac{p_1}{p_2}, \ z_6=q_2+\frac{q_1p_1+\beta_6-(\beta_0+\beta_4)}{p_2}-\frac{s-1}{p_2^2}, \ w_6=p_2,
\end{align*}
where $r_2^{04} \left(H_{04}^{(1)}-\frac{1}{p_1} \right), \ r_2^{04} \left(H_{04}^{(2)}+\frac{1}{p_1} \right), \ r_3^{04} \left(H_{04}^{(2)}+\frac{1}{p_2} \right), \ r_6^{04} \left(H_{04}^{(2)}+\frac{1}{p_2} \right)$.

\hspace{0.5cm}

{\bf Hamiltonians $H_{045}^{(1)}=r_5^{04}(H_{04}^{(1)}), \ H_{045}^{(2)}=r_5^{04}(H_{04}^{(2)})$}

\begin{align*}
&r_0^{045}:x_0=\frac{1}{q_1}, \ y_0=-(p_1 q_1+\beta_0+\beta_5)q_1, \ z_0=q_2, \ w_0=p_2, \\
&r_1^{045}:x_1=\frac{1}{q_1}, \ y_1=-(p_1 q_1+\beta_1+\beta_0+\beta_5)q_1, \ z_1=q_2, \ w_1=p_2, \\
&r_2^{045}:x_2=x_5 w_5^2, \ y_2=\frac{y_5}{w_5^2}, \ z_2=z_5+\frac{2x_5 y_5+\beta_2-(\beta_0+\beta_4)}{w_5}+\frac{t-s}{w_5^2} \ w_2=w_5, \\
&r_3^{045}:x_3=\frac{q_1}{p_2}, \ y_3=p_1 p_2, \ z_3=q_2-\frac{q_1p_1-\beta_3+\beta_0+\beta_4+\beta_5}{p_2}-\frac{s}{p_2^2}, \ w_3=p_2, \\
&r_4^{045}:x_4=q_1, \ y_4=p_1-\frac{\beta_4 q_2}{q_1q_2-1}, \ z_4=q_2, \ w_4=p_2-\frac{\beta_4 q_1}{q_1q_2-1}, \\
&r_5^{045}:x_5=\frac{1}{q_1}, \ y_5=-(p_1 q_1+\beta_5)q_1, \ z_5=q_2, \ w_5=p_2,\\
&r_6^{045}:x_6=\frac{q_1}{p_2}, \ y_6=p_1 p_2, \ z_6=q_2-\frac{q_1p_1-\beta_6+\beta_0+\beta_4+\beta_5}{p_2}-\frac{s-1}{p_2^2}, \ w_6=p_2,
\end{align*}
where $r_2^{045} \left(r_5^{045}(H_{045}^{(1)})-\frac{1}{w_5} \right), \ r_2^{045} \left(r_5^{045}(H_{045}^{(2)})+\frac{1}{w_5} \right), \ r_3^{045} \left(H_{045}^{(2)}+\frac{1}{p_2} \right), \ r_6^{045} \left(H_{045}^{(2)}+\frac{1}{p_2} \right)$.

{\bf Hamiltonians $H_{15}^{(1)}=r_{01}\left(-\frac{1}{T^2} H_1^{FS} \right), \ H_{15}^{(2)}=r_{15}\left(-\frac{1}{S^2} H_2^{FS} \right), \ r_{15}:Q_1=\frac{1-q_1}{q_1}, \ P_1=-(p_1 q_1+\beta_1)q_1, \ Q_2=\frac{1-q_2}{q_2}, \ P_2=-(p_2 q_2+\beta_5)q_2, \ T=\frac{1-t}{t}, \ S=\frac{1-s}{s}$}
\begin{align*}
&r_0^{15}:x_0=-((q_1-q_2)p_1-\beta_0)p_1, \ y_0=\frac{1}{p_1}, \ z_0=q_2, \ w_0=p_2+p_1, \\
&r_1^{15}:x_1=\frac{1}{q_1}, \ y_1=-(q_1p_1+\beta_1)q_1, \ z_1=q_2, \ w_1=p_2, \\
&r_2^{15}:x_2=-((q_1-t)p_1-\beta_2)p_1, \ y_2=\frac{1}{p_1}, \ z_2=q_2, \ w_2=p_2, \\
&r_3^{15}:x_3=\frac{1}{q_1}, \ y_3=-(p_1 q_1+p_2 q_2+\beta_3+\beta_1+\beta_5)q_1, \ z_3=\frac{q_2}{q_1}, \ w_3=p_2 q_1, \\
&r_4^{15}:x_4=q_1, \ y_4=p_1, \ z_4=-((q_2-s)p_2-\beta_4)p_2, \ w_4=\frac{1}{p_2}, \\
&r_5^{15}:x_5=q_1, \ y_5=p_1, \ z_5=\frac{1}{q_2}, \ w_5=-(p_2q_2+\beta_5)q_2,\\
&r_6^{15}:x_6=-(q_1 p_1+q_2 p_2-\beta_6)p_1, \ y_6=\frac{1}{p_1}, \ z_6=q_2 p_1, \ w_6=\frac{p_2}{p_1},\\
\end{align*}
where $r_2^{15} \left(H_{15}^{(1)}-p_1 \right), \ r_4^{15} \left(H_{15}^{(2)}-p_2 \right)$.

{\bf Hamiltonians $H_{150}^{(1)}=r_0^{15}\left(H_{15}^{(1)} \right), \ H_{150}^{(2)}=r_0^{15}\left(H_{15}^{(2)} \right)$}
\begin{align*}
&r_0^{150}:x_0=-(q_1 p_1-\beta_0)p_1, \ y_0=\frac{1}{p_1}, \ z_0=q_2, \ w_0=p_2, \\
&r_1^{150}:x_1=-(q_1 p_1-\beta_1-\beta_0)p_1, \ y_1=\frac{1}{p_1}, \ z_1=q_2, \ w_1=p_2, \\
&r_2^{150}:x_2=q_1+\frac{\beta_2-\beta_0}{p_1}+\frac{t-q_2}{p_1^2}, \ y_2=p_1, \ z_2=q_2, \ w_2=p_2-\frac{1}{p_1}, \\
&r_3^{150}:x_3=-(q_1 p_1-q_2 p_2-\beta_3-\beta_0-\beta_1-\beta_5)p_1, \ y_3=\frac{1}{p_1}, \ z_3=\frac{q_2}{p_1}, \ w_3=p_2 p_1, \\
&r_4^{150}:x_4=q_1, \ y_4=p_1, \ z_4=-((q_2-s)p_2-\beta_4)p_2, \ w_4=\frac{1}{p_2}, \\
&r_5^{150}:x_5=q_1+\frac{\beta_5 p_2}{p_1 p_2-1}, \ y_5=p_1, \ z_5=q_2+\frac{\beta_5 p_1}{p_1 p_2-1}, \ w_5=p_2,\\
&r_6^{150}:x_6=\frac{1}{q_1}, \ y_6=-(p_1 q_1-p_2 q_2+\beta_6-\beta_0)q_1, \ z_6=q_2 q_1, \ w_6=\frac{p_2}{q_1},\\
\end{align*}
where $r_2^{150} \left(H_{150}^{(1)}-\frac{1}{p_1} \right), \ r_4^{150} \left(H_{150}^{(2)}-p_2 \right)$.

{\bf Hamiltonians $H_{1503}^{(1)}=sr_3 \left(H_{150}^{(1)} \right), \ H_{1503}^{(2)}=sr_3 \left(H_{150}^{(2)} \right), \ sr_3:Q_1=q_1, \ P_1=p_1-\frac{q_2 p_2+\beta_3+\beta_0+\beta_1+\beta_5}{q_1}, \ Q_2=q_2 q_1, \ P_2=\frac{p_2}{q_1}$}
\begin{align*}
&r_0^{1503}:x_0=-(q_1p_1+q_2p_2-\delta_0)p_1, \ y_0=\frac{1}{p_1}, \ z_0=q_2p_1, \ w_0=\frac{p_2}{p_1}, \\
&r_1^{1503}:x_1=-(q_1p_1+q_2p_2-\delta_1-\delta_0)p_1, \ y_1=\frac{1}{p_1}, \ z_1=q_2p_1, \ w_1=\frac{p_2}{p_1}, \\
&r_2^{1503}:x_2=q_1+\frac{q_2}{t}+\frac{2q_2p_2+2\delta_2}{p_1}+\frac{t}{p_1^2}, \ y_2=p_1, \ z_2=q_2p_1^2, \ w_2=\frac{p_2+\frac{p_1}{t}}{p_1^2}, \\
&r_3^{1503}:x_3=\frac{1}{q_1}, \ y_3=-(p_1q_1+\delta_3)q_1, \ z_3=q_2, \ w_3=p_2, \\
&r_4^{1503}:x_4=-\left( \left(q_1-\frac{1}{s} q_2 \right)p_1-\delta_4 \right)p_1, \ y_4=\frac{1}{p_1}, \ z_4=q_2, \ w_4=p_2+\frac{1}{s} p_1, \\
&r_5^{1503}:x_5=q_1p_2, \ y_5=\frac{p_1}{p_2}, \ z_5=q_2+\frac{q_1p_1+2\delta_5}{p_2}-\frac{1}{p_2^2}, \ w_5=p_2,\\
&r_6^{1503}:x_6=-(q_1p_1+q_2p_2-\delta_6-\delta_1-\delta_0)p_1, \ y_6=\frac{1}{p_1}, \ z_6=q_2 q_1, \ w_6=\frac{p_2}{q_1},\\
\end{align*}
where $r_2^{1503} \left(H_{1503}^{(1)}+\frac{p_1 q_2}{t^2}-\frac{1}{p_1} \right), \ r_4^{1503} \left(H_{1503}^{(2)}+\frac{p_1 q_2}{s^2} \right)$.

Here, $\delta_j \ (j=0,1,\ldots,6)$ are complex constant parameters satisfying the parameter's relation\rm{: \rm}
\begin{align}
\begin{split}
&\beta_0=-\delta_6, \ \beta_1=-\delta_1, \ \beta_2=2\delta_0+2\delta_1+2\delta_2+\delta_6,\\
&\beta_3=-2\delta_0-\delta_1-2\delta_5-\delta_6, \ \beta_4=\delta_4, \ \beta_5=\delta_0+\delta_1+2\delta_5+\delta_6, \ \beta_6=\delta_0+\delta_1+\delta_3,
\end{split}
\end{align}
\begin{align}
\begin{split}
&3\delta_0+2\delta_1+2\delta_2+\delta_3+\delta_4+2\delta_5+\delta_6=1.
\end{split}
\end{align}
We remark that the transformations $r_2^{1503},r_5^{1503}$ are not its auto-B{\"a}cklund transformations. It is still an open question whether the transformations $r_2^{1503},r_5^{1503}$ can be considered as each transformation denoted by the symbol $\odot$ in the Oshima's paper (see \cite{Oshima}).

It is also still an open question whether we can obtain the Hamiltonian system with $H_{1503}^{(1)},H_{1503}^{(2)}$ by solving $3 \times 3$ Lax pair  (cf. \cite{KNS,Oshima}) satisfying the following Riemann scheme:
\begin{equation*}
\begin{pmatrix}
X=0 & X=1 &  X=t &  X=\infty\\
\begin{matrix}
0 \\
\theta_1^{0} \\
\theta_2^{0}
\end{matrix}  & \begin{matrix}
0 \\
0 \\
\theta^{1}
\end{matrix}  &  \begin{matrix}
0 \\
0 \\
\theta^{t}
\end{matrix} &  \begin{matrix}
\theta_1^{\infty}\\
\theta_2^{\infty} \\
\theta_3^{\infty}
\end{matrix}
\end{pmatrix}.
\end{equation*}

\section{Appendix E: 5-parameter family of polynomial Hamiltonian system}

{\bf Holomorphy conditions}

Define birational and symplectic transformations $r_i \ (i=0,1,\ldots,5)$ as follows:

\begin{align}\label{holo;5p}
\begin{split}
r_0:(x_0,y_0,z_0,w_0)=&\left(\frac{1}{q_1},-(q_1p_1+\alpha_0)q_1,q_2,p_2 \right),\\
r_1:(x_1,y_1,z_1,w_1)=&\left(q_1,p_1,\frac{1}{q_2},-(p_2q_2+\alpha_1)q_2 \right),\\
r_2:(x_2,y_2,z_2,w_2)=&\left(-((q_1-q_2)p_1-\alpha_2)p_1,\frac{1}{p_1},q_2,p_2+p_1 \right),\\
r_3:(x_3,y_3,z_3,w_3)=&\left(-((q_1-1)p_1+(q_2-1)p_2-\alpha_3)p_1,\frac{1}{p_1},(q_2-1)p_1,\frac{p_2}{p_1} \right),\\
r_4:(x_4,y_4,z_4,w_4)=&\left(-(q_1p_1+q_2p_2-\alpha_4)p_1,\frac{1}{p_1},q_2p_1,\frac{p_2}{p_1} \right),\\
r_5:(x_5,y_5,z_5,w_5)=&\left(-((q_1-t)p_1+(q_2-s)p_2-\alpha_5)p_1,\frac{1}{p_1},(q_2-s)p_1,\frac{p_2}{p_1} \right).
\end{split}
\end{align}
There exist two polynomials $H_1$ and $H_2$, such that the Hamiltonian system
\begin{equation}\label{eq:5p}
  \left\{
  \begin{aligned}
   dq_1 =&\frac{\partial H_1}{\partial p_1}dt+\frac{\partial H_2}{\partial p_1}ds, \quad dp_1 =-\frac{\partial H_1}{\partial q_1}dt-\frac{\partial H_2}{\partial q_1}ds,\\
   dq_2 =&\frac{\partial H_1}{\partial p_2}dt+\frac{\partial H_2}{\partial p_2}ds, \quad dp_2 =-\frac{\partial H_1}{\partial q_2}dt-\frac{\partial H_2}{\partial q_2}ds
   \end{aligned}
  \right.
\end{equation}
is transformed into a polynomial Hamiltonian system under the action of each $r_i \ (i=0,1,\ldots,5)$, where two polynomial Hamiltonians $H_1,H_2$ are given by {\rm (cf. \cite{FS1,KNS,Fuji}) \rm}

{\footnotesize
\begin{align}
\begin{split}
&(3\alpha_0+3\alpha_1+2\alpha_2+\alpha_3+\alpha_4+2\alpha_5)t(t-1)(t-s)H_1\\
&=\tilde{H}_{VI}(q_1,p_1,t,s;\alpha_5,\alpha_0+3\alpha_1+2\alpha_2+\alpha_5,\alpha_0,\alpha_3+\alpha_1,\alpha_4-\alpha_1)\\
&+\tilde{H}_{VI}(q_2,p_2,t,s;3\alpha_0+\alpha_1+2\alpha_2+\alpha_5,\alpha_5,\alpha_1,\alpha_3+\alpha_0,\alpha_4-\alpha_0)\\
&+\{ 2 q_1^2 q_2^2 - ((3 t - s + 1) q_1 - t (s + 2) + 2 s t ) q_2^2 - ((3 s - t + 1) q_2 - s (t + 2) + 2 t s) q_1^2 \\
&+ s t (4 q_1 q_2 - q_1 - q_2) \} p_1 p_2 + \alpha_0 \{ 2 q_1 q_2^2 + (-3 s + t - 1) q_1 q_2 - (t - 2) s q_1 + (s - 2) t q_2 + s t \} p_2\\
& + \alpha_1 \{ 2 q_1^2 q_2 + (-3 t + s - 1) q_1 q_2 + s (t - 2) q_1 - t (s - 2) q_2 + s t \} p_1 + 2 \alpha_0 \alpha_1 q_1 q_2,\\
H_2 &=\pi(H_1), \quad \pi=\{q_1 \leftrightarrow q_2, \ p_1 \leftrightarrow p_2, \  t \leftrightarrow s, \alpha_0 \leftrightarrow \alpha_1 \}.
\end{split}
\end{align}
The symbol $\tilde{H}_{VI}(q,p,t,\eta;\gamma_0,\gamma_1,\gamma_2,\gamma_3,\gamma_4)$ denotes (see \cite{sasa3})
\begin{align}
\begin{split}
&\tilde{H}_{VI}(q,p,t,\eta;\gamma_0,\gamma_1,\gamma_2,\gamma_3,\gamma_4)\\
&=q(q-1)(q-\eta)(q-t)p^2+\{\gamma_1(t-\eta)q(q-1)+2\gamma_2 q(q-1)(q-\eta)\\
&+\gamma_3 (t-1)q(q-\eta)+\gamma_4 t(q-1)(q-\eta) \}p\\
&+\gamma_2 \{(\gamma_1+\gamma_2)(t-\eta)+\gamma_2(q-1)+\gamma_3(t-1)+t \gamma_4 \}q, \quad (\gamma_0+\gamma_1+2\gamma_2+\gamma_3+\gamma_4=1).
\end{split}
\end{align}
}

\vspace{0.5cm}

We note that the holomorphy conditions should be read that in the Hamiltonian $H_1$
\begin{align*}
\begin{split}
&r_5(H_1 - p_1)
\end{split}
\end{align*}
are polynomials with respect to $x_5,y_5,z_5,w_5$, and in the Hamiltonian $H_2$
\begin{align*}
\begin{split}
&r_5(H_2 -p_2)
\end{split}
\end{align*}
are polynomials with respect to $x_5,y_5,z_5,w_5$.

\vspace{0.5cm}

{\bf Completely integrable}

{\footnotesize
\begin{proposition}
Setting
\begin{equation}
K_1:=-H_1+\frac{\alpha_4(\alpha_0+\alpha_1)( \rm{Log \rm}(s-t)-\rm{Log \rm}(s-1) )}{(3\alpha_0+3\alpha_1+2\alpha_2+\alpha_3+\alpha_4+2\alpha_5)(t-1)^2}, \quad K_2:=-H_2.
\end{equation}
Two Hamiltonians $K_1$ and $K_2$ satisfy
\begin{equation}
\{K_1,K_2\}+\left(\frac{\partial}{\partial s} \right)K_1-\left(\frac{\partial}{\partial t} \right)K_2=0,
\end{equation}
where  $\{,\}$ denotes the Poisson brackets:
\begin{equation}
\{L_1,L_2\}=\frac{\partial L_1}{\partial p_1}\frac{\partial L_2}{\partial q_1}-\frac{\partial L_1}{\partial q_1}\frac{\partial L_2}{\partial p_1}+\frac{\partial L_1}{\partial p_2}\frac{\partial L_2}{\partial q_2}-\frac{\partial L_1}{\partial q_2}\frac{\partial L_2}{\partial p_2}.
\end{equation}
\end{proposition}
}

{\bf Symmetry}

With {\it the notation} $(*):=(q_1,p_1,q_2,p_2,t,s;\alpha_0,\alpha_1,\alpha_2,\alpha_3,\alpha_4,\alpha_5)$;

{\tiny
\begin{align}
\begin{split}
&s_0:(*) \rightarrow \left(q_1+\frac{\alpha_0}{p_1},p_1,q_2,p_2,t,s;-\alpha_0,\alpha_1,\alpha_2+\alpha_0,\alpha_3+\alpha_0,\alpha_4+\alpha_0,\alpha_5+\alpha_0 \right),\\
&s_1:(*) \rightarrow \left(q_1,p_1,q_2+\frac{\alpha_1}{p_2},p_2,t,s;\alpha_0,-\alpha_1,\alpha_2+\alpha_1,\alpha_3+\alpha_1,\alpha_4+\alpha_1,\alpha_5+\alpha_1 \right),\\
&s_2:(*) \rightarrow \left(q_1,p_1-\frac{\alpha_2}{q_1-q_2},q_2,p_2+\frac{\alpha_2}{q_1-q_2},t,s;\alpha_0+\alpha_2,\alpha_1+\alpha_2,-\alpha_2,\alpha_3-\alpha_2,\alpha_4-\alpha_2,\alpha_5 \right),\\
&s_3:(*) \rightarrow \left(1-q_1,-p_1,1-q_2,-p_2,1-t,1-s;\alpha_0,\alpha_1,\alpha_2,\alpha_4,\alpha_3,\alpha_5 \right),\\
&s_4:(*) \rightarrow (\frac{q_1 (-q_1 p_1 + q_1^2 p_1 - p_2 q_2 + p_2 q_2^2 - \alpha_0 + 
     q_1 \alpha_0 - \alpha_1 + q_2 \alpha_1 - \alpha_3)}{g_1},\\
& \frac{-g_1 (q_1 p_1^2 - q_1^2 p_1^2 + 
        p_2 p_1 q_2 - p_2 p_1 q_2^2 + \alpha_0^2 + q_1 p_1 \alpha_1 - 
        p_1 q_2 \alpha_1 + \alpha_0 \alpha_1 + 
        q_1 p_1 \alpha_3 + \alpha_0 \alpha_3 - p_1 \alpha_4 + 
        q_1 p_1 \alpha_4 + \alpha_0 \alpha_4)}{(-q_1 p_1 + q_1^2 p_1 - p_2 q_2 + 
       p_2 q_2^2 - \alpha_0 + q_1 \alpha_0 - \alpha_1 + 
       q_2 \alpha_1 - \alpha_3) (-q_1 p_1 + q_1^2 p_1 - p_2 q_2 + p_2 q_2^2 + 
       q_1 \alpha_0 + q_2 \alpha_1 + \alpha_4)}, \\
& \frac{q_2 (-q_1 p_1 + q_1^2 p_1 - p_2 q_2 + p_2 q_2^2 - \alpha_0 + 
     q_1 \alpha_0 - \alpha_1 + q_2 \alpha_1 - \alpha_3)}{g_2 }, \\
& \frac{ -g_2 (p_2 q_1 p_1 - p_2 q_1^2 p_1 + 
        p_2^2 q_2 - p_2^2 q_2^2 - p_2 q_1 \alpha_0 + 
        p_2 q_2 \alpha_0 + \alpha_0 \alpha_1 + \alpha_1^2 + 
        p_2 q_2 \alpha_3 + \alpha_1 \alpha_3 - p_2 \alpha_4 + 
        p_2 q_2 \alpha_4 + \alpha_1 \alpha_4)}{(-q_1 p_1 + q_1^2 p_1 - p_2 q_2 + 
       p_2 q_2^2 - \alpha_0 + q_1 \alpha_0 - \alpha_1 + q_2 \alpha_1 - \alpha_3) (-q_1 p_1 + q_1^2 p_1 - p_2 q_2 + p_2 q_2^2 + q_1 \alpha_0 + q_2 \alpha_1 + \alpha_4)},\\
&t,s;\alpha_0,\alpha_1,\alpha_2,-\alpha_0-\alpha_1-\alpha_4,-\alpha_0-\alpha_1-\alpha_3,\alpha_0+\alpha_1+\alpha_3+\alpha_4+\alpha_5 ),\\
&s_5:(*) \rightarrow (\frac{q_1 (-s t q_1 p_1 + s q_1^2 p_1 - s t p_2 q_2 + t p_2 q_2^2 - s t \alpha_0 + 
     s q_1 \alpha_0 - s t \alpha_1 + t q_2 \alpha_1 - 
     s t \alpha_5)}{g_3}, \\
& \frac{-g_3 (s t q_1 p_1^2 - s q_1^2 p_1^2 + s t p_2 p_1 q_2 - 
        t p_2 p_1 q_2^2 + s \alpha_0^2 + s q_1 p_1 \alpha_1 - 
        t p_1 q_2 \alpha_1 + s \alpha_0 \alpha_1 - s t p_1 \alpha_4 + 
        s q_1 p_1 \alpha_4 + s \alpha_0 \alpha_4 + s q_1 p_1 \alpha_5 + 
        s \alpha_0 \alpha_5)}{(-s t q_1 p_1 + s q_1^2 p_1 - s t p_2 q_2 + 
       t p_2 q_2^2 + s q_1 \alpha_0 + t q_2 \alpha_1 + 
       s t \alpha_4) (-s t q_1 p_1 + s q_1^2 p_1 - s t p_2 q_2 + t p_2 q_2^2 - 
       s t \alpha_0 + s q_1 \alpha_0 - s t \alpha_1 + t q_2 \alpha_1 -
        s t \alpha_5)}, \\
& \frac{q_2 (-s t q_1 p_1 + s q_1^2 p_1 - s t p_2 q_2 + t p_2 q_2^2 - s t \alpha_0 + 
     s q_1 \alpha_0 - s t \alpha_1 + t q_2 \alpha_1 - 
     s t \alpha_5)}{g_4}, \\
& \frac{-g_4 (s t p_2 q_1 p_1 - s p_2 q_1^2 p_1 + s t p_2^2 q_2 - 
        t p_2^2 q_2^2 - s p_2 q_1 \alpha_0 + t p_2 q_2 \alpha_0 + 
        t \alpha_0 \alpha_1 + t \alpha_1^2 - s t p_2 \alpha_4 + 
        t p_2 q_2 \alpha_4 + t \alpha_1 \alpha_4 + t p_2 q_2 \alpha_5 + 
        t \alpha_1 \alpha_5)}{(-s t q_1 p_1 + s q_1^2 p_1 - s t p_2 q_2 + 
       t p_2 q_2^2 + s q_1 \alpha_0 + t q_2 \alpha_1 + 
       s t \alpha_4) (-s t q_1 p_1 + s q_1^2 p_1 - s t p_2 q_2 + t p_2 q_2^2 - 
       s t \alpha_0 + s q_1 \alpha_0 - s t \alpha_1 + t q_2 \alpha_1 -
        s t \alpha_5)},\\
&t,s;\alpha_0,\alpha_1,\alpha_0+\alpha_1+\alpha_2+\alpha_4+\alpha_5,\alpha_0+\alpha_1+\alpha_3+\alpha_4+\alpha_5,-\alpha_0-\alpha_1-\alpha_5,-\alpha_0-\alpha_1-\alpha_4 ),\\
&s_6:(*) \rightarrow \left(q_2,p_2,q_1,p_1,s,t;\alpha_1,\alpha_0,\alpha_2,\alpha_3,\alpha_4,\alpha_5 \right),
\end{split}
\end{align}}
where 
\begin{align*}
\begin{split}
&g_1:=-q_1 p_1 + q_1^2 p_1 - p_2 q_2 + p_2 q_2^2 - q_1 \alpha_1 + q_2 \alpha_1 - q_1 \alpha_3 + \alpha_4 - q_1 \alpha_4,\\
&g_2:=-q_1 p_1 + q_1^2 p_1 - p_2 q_2 + p_2 q_2^2 + q_1 \alpha_0 - q_2 \alpha_0 - q_2 \alpha_3 + \alpha_4 - q_2 \alpha_4,\\
&g_3:=-s t q_1 p_1 + s q_1^2 p_1 - s t p_2 q_2 + t p_2 q_2^2 - s q_1 \alpha_1 + t q_2 \alpha_1 + s t \alpha_4 - s q_1 \alpha_4 - s q_1 \alpha_5,\\
&g_4:=-s t q_1 p_1 + s q_1^2 p_1 - s t p_2 q_2 + t p_2 q_2^2 + s q_1 \alpha_0 - t q_2 \alpha_0 + s t \alpha_4 - t q_2 \alpha_4 - t q_2 \alpha_5.
\end{split}
\end{align*}

Finally, we remark some holomorphy conditions of the system \eqref{eq:5p}.

{\bf Hamiltonians $L_{01}^{(1)}=r_{01}\left(-\frac{1}{T^2} H_1 \right), \ L_{01}^{(2)}=r_{01}\left(-\frac{1}{S^2} H_2 \right), \ r_{01}:Q_1=\frac{1}{q_1}, \ P_1=-(q_1p_1+\alpha_0)q_1, \ Q_2=\frac{1}{q_2}, \ P_2=-(p_2q_2+\alpha_1)q_2, \ T=\frac{1}{t}, \ S=\frac{1}{s}$}
\begin{align*}
&r_0^{01}:x_0=\frac{1}{q_1}, \ y_0=-(q_1p_1+\alpha_0)q_1, \ z_0=q_2, \ w_0=p_2, \\
&r_1^{01}:x_1=q_1, \ y_1=p_1, \ z_1=\frac{1}{q_2}, \ w_1=-(p_2q_2+\alpha_1)q_2, \\
&r_2^{01}:x_2=-((q_1-q_2)p_1-\alpha_2)p_1, \ y_2=\frac{1}{p_1}, \ z_2=q_2, \ w_2=p_2+p_1, \\
&r_3^{01}:x_3=-((q_1-1)p_1+(q_2-1)p_2-\alpha_3)p_1, \ y_3=\frac{1}{p_1}, \ z_3=(q_2-1)p_1, \ w_3=\frac{p_2}{p_1}, \\
&r_4^{01}:x_4=\frac{1}{q_1}, \ y_4=-(p_1 q_1+p_2 q_2+\alpha_4+\alpha_0+\alpha_1)q_1, \ z_4=\frac{q_2}{q_1}, \ w_4=p_2 q_1, \\
&r_5^{01}:x_5=-((q_1-t)p_1+(q_2-s)p_2-\alpha_5)p_1, \ y_5=\frac{1}{p_1}, \ z_5=(q_2-s)p_1, \ w_5=\frac{p_2}{p_1},
\end{align*}
where $r_5^{01} \left(L_{01}^{(1)}-p_1 \right), \ r_5^{01} \left(L_{01}^{(2)}-p_2 \right)$.

Here, for notational convenience, we have renamed $(Q_1,P_1,Q_2,P_2,T,S)$ to $(q_1,p_1,q_2,p_2,t,s)$ (which are not the same as the previous $(q_1,p_1,q_2,p_2,t,s)$).

{\bf Hamiltonians $L_{014}^{(1)}=r_4^{01}\left(-\frac{1}{T^2} L_{01}^{(1)}-\frac{S}{T^2} L_{01}^{(2)} \right), \ L_{014}^{(2)}=r_4^{01}\left(\frac{1}{T} L_{01}^{(2)} \right), \ r_4^{01}:Q_1=\frac{1}{q_1}, \ P_1=-(p_1 q_1+p_2 q_2+\alpha_4+\alpha_0+\alpha_1)q_1, \ Q_2=\frac{q_2}{q_1}, \ P_2=p_2 q_1, \ T=\frac{1}{t}, \ S=\frac{s}{t}$}
\begin{align*}
&r_0^{014}:x_0=q_1 p_2, \ y_0=\frac{p_1}{p_2}, \ z_0=-(p_2 q_2+p_1 q_1+\alpha_1+\alpha_4)p_2, \ w_0=\frac{1}{p_2}, \\
&r_1^{014}:x_1=q_1, \ y_1=p_1, \ z_1=\frac{1}{q_2}, \ w_1=-(p_2q_2+\alpha_1)q_2, \\
&r_2^{014}:x_2=q_1, \ y_2=p_1, \ z_2=-((q_2-1)p_2-\alpha_2)p_2, \ w_2=\frac{1}{p_2}, \\
&r_3^{014}:x_3=-((q_1-1)p_1+(q_2-1)p_2-\alpha_3)p_1, \ y_3=\frac{1}{p_1}, \ z_3=(q_2-1)p_1, \ w_3=\frac{p_2}{p_1}, \\
&r_4^{014}:x_4=\frac{1}{q_1}, \ y_4=-(p_1 q_1+p_2 q_2+\alpha_4+\alpha_0+\alpha_1)q_1, \ z_4=\frac{q_2}{q_1}, \ w_4=p_2 q_1, \\
&r_5^{014}:x_5=-((q_1-t)p_1+(q_2-s)p_2-\alpha_5)p_1, \ y_5=\frac{1}{p_1}, \ z_5=(q_2-s)p_1, \ w_5=\frac{p_2}{p_1},
\end{align*}
where $r_5^{014} \left(L_{014}^{(1)}-p_1 \right), \ r_5^{014} \left(L_{014}^{(2)}-p_2 \right)$.

{\bf Hamiltonians $L_1=r_0^{014}\left(L_{014}^{(1)} \right), \ L_2=r_0^{014}\left(L_{014}^{(2)} \right)$}
\begin{align*}
&R_0:x_0=q_1 p_2, \ y_0=\frac{p_1}{p_2}, \ z_0=-(p_2 q_2+p_1 q_1+\alpha_1+\alpha_4)p_2, \ w_0=\frac{1}{p_2}, \\
&R_1:x_1=q_1, \ y_1=p_1, \ z_1=-(q_2 p_2-\alpha_0)p_2, \ w_1=\frac{1}{p_2}, \\
&R_2:x_2=q_1 p_2, \ y_2=\frac{p_1}{p_2}, \ z_2=q_2+\frac{q_1 p_1+\alpha_2+\alpha_1+\alpha_4}{p_2}+\frac{1}{p_2^2}, \ w_2=p_2, \\
&R_3:x_3=q_1-\frac{1}{p_2}, \ y_3=p_1, \ z_3=q_2+\frac{\alpha_3+\alpha_1+\alpha_4}{p_2}+\frac{p_1+1}{p_2^2}, \ w_3=p_2, \\
&R_4:x_4=q_1 p_2, \ y_4=\frac{p_1}{p_2}, \ z_4=-(p_2 q_2+p_1 q_1+\alpha_4)p_2, \ w_4=\frac{1}{p_2}, \\
&R_5:x_5=q_1-\frac{t}{p_2}, \ y_5=p_1, \ z_5=q_2+\frac{\alpha_5+\alpha_1+\alpha_4}{p_2}+\frac{t p_1+s}{p_2^2}, \ w_5=p_2,
\end{align*}
where $R_5 \left(L_1-\frac{p_1}{p_2} \right), \ R_5 \left(L_2-\frac{1}{p_2} \right)$. We remark that the transformation $R_2$ is not its auto-B{\"a}cklund transformation. It is still an open question whether the transformation $R_2$ can be considered as a transformation denoted by the symbol $\odot$ in the Oshima's paper (see \cite{Oshima}).

{\scriptsize
It is still an open question whether we can obtain the Hamiltonian system with $L_1,L_2$ by solving $4 \times 4$ Lax pair  (cf. \cite{KNS,Oshima}) satisfying the following Riemann scheme:
\begin{equation}
\begin{pmatrix}
X=0 & X=1 &  X=t &  X=\infty\\
\begin{matrix}
0 \\
0 \\
\theta_1 \\
\theta_1
\end{matrix}  & \begin{matrix}
0 \\
0 \\
\theta_2 \\
\theta_2
\end{matrix}  &  \begin{matrix}
0 \\
0 \\
\theta_3 \\
\theta_3
\end{matrix} &  \begin{matrix}
0\\
\alpha_4 \\
\alpha_4+\alpha_1 \\
\alpha_0
\end{matrix}
\end{pmatrix}.
\end{equation}

Here, we will conjecture the following relations between Riemann data and Holomorphy conditions $R_i \ (i=0,1,\ldots,5);$

$\begin{pmatrix}
X=0,1,t \\
\begin{matrix}
0 \\
0 \\
\theta \\
\theta
\end{matrix}
\end{pmatrix} \Longleftrightarrow $ Holomorphy conditions $R_2,R_3,R_5$, \ $\begin{pmatrix}
X=\infty \\
 \begin{matrix}
0\\
\alpha_4 \\
\alpha_4+\alpha_1 \\
\alpha_0
\end{matrix}
\end{pmatrix} \Longleftrightarrow $ Holomorphy conditions $\begin{pmatrix}
\begin{matrix}
R_4 \\
R_0 \\
R_1
\end{matrix}
\end{pmatrix}.$

Setting $(X,Y,Z,W)=(q_1 p_2,p_1,q_2 p_2^2,p_2)$, we will see the following relations;

Holomorphy condition $R_2$ $\Longleftrightarrow$ Accessible sing. $(X,Y,Z,W)=(X,0,-1,0)$, \ Local index:$\begin{pmatrix}
\begin{matrix}
0 \\
1 \\
2 \\
1
\end{matrix}
\end{pmatrix}$, 

Holomorphy condition $R_3$ $\Longleftrightarrow$ Accessible sing. $(X,Y,Z,W)=(1,Y,-1-Y,0)$, \ Local index:$\begin{pmatrix}
\begin{matrix}
1 \\
0 \\
2 \\
1
\end{matrix}
\end{pmatrix}$,

Holomorphy condition $R_5$ $\Longleftrightarrow$ Accessible sing. $(X,Y,Z,W)=(t,Y,-s-tY,0)$, \ Local index:$\begin{pmatrix}
\begin{matrix}
1 \\
0 \\
2 \\
1
\end{matrix}
\end{pmatrix}$.

We remark that

well-known Holomorphy condition $R:x=q_1, \ y=p_1, \ z=q_2+\frac{\alpha}{p_2}+\frac{\eta}{p_2^2}, \ w=p_2$

$\Longleftrightarrow$ Accessible sing. $(X,Y,Z,W)=(0,Y,-\eta,0)$, \ Local index:$\begin{pmatrix}
\begin{matrix}
1 \\
0 \\
2 \\
1
\end{matrix}
\end{pmatrix}$.

We note that we will make Holomorphy condition $R_2$.

{\bf Step 0:} We make a change of variables:
$$
X=q_1 p_2, \quad Y=p_1, \quad Z=q_2 p_2^2, \quad W=p_2.
$$

{\bf Step 1:} We blow up along the curve $\{(X,Y,Z,W)|Y=0, \ Z=-1, \ W=0\}$:
$$
X^{(1)}=X, \quad Y^{(1)}=\frac{Y}{W}, \quad Z^{(1)}=\frac{Z+1}{W}, \quad W^{(1)}=W.
$$

{\bf Step 2:} We blow up along the surface $\{(X^{(1)},Y^{(1)},Z^{(1)},W^{(1)})|Z^{(1)}=-(X^{(1)} Y^{(1)}+\alpha_2+\alpha_1+\alpha_4), \  W^{(1)}=0\}$:
$$
X^{(2)}=X^{(1)}, \quad Y^{(2)}=Y^{(1)}, \quad Z^{(2)}=\frac{Z^{(1)}+X^{(1)} Y^{(1)}+\alpha_2+\alpha_1+\alpha_4}{W^{(1)}}, \quad W^{(2)}=W^{(1)}.
$$
By taking the coordinate system as
$$
(x_2,y_2,z_2,w_2)=(X^{(2)},Y^{(2)},Z^{(2)},W^{(2)}),
$$
we can obtain the coordinate system $R_2$.

We also note that we will make Holomorphy condition $R_3$.

{\bf Step 0:} We make a change of variables:
$$
X=q_1 p_2, \quad Y=p_1, \quad Z=q_2 p_2^2, \quad W=p_2.
$$

{\bf Step 1:} We blow up along the curve $\{(X,Y,Z,W)|X=1, \ Z=-1-Y, \ W=0\}$:
$$
X^{(1)}=\frac{X-1}{W}, \quad Y^{(1)}=Y, \quad Z^{(1)}=\frac{Z+1+Y}{W}, \quad W^{(1)}=W.
$$

{\bf Step 2:} We blow up along the surface $\{(X^{(1)},Y^{(1)},Z^{(1)},W^{(1)})|Z^{(1)}=-(\alpha_3+\alpha_1+\alpha_4), \  W^{(1)}=0\}$:
$$
X^{(2)}=X^{(1)}, \quad Y^{(2)}=Y^{(1)}, \quad Z^{(2)}=\frac{Z^{(1)}+\alpha_3+\alpha_1+\alpha_4}{W^{(1)}}, \quad W^{(2)}=W^{(1)}.
$$
By taking the coordinate system as
$$
(x_3,y_3,z_3,w_3)=(X^{(2)},Y^{(2)},Z^{(2)},W^{(2)}),
$$
we can obtain the coordinate system $R_3$.

It is still an open question whether in Lax pair theory we can distinguish Holomorphy conditions $R,R_3,R_5$ from Holomorphy condition $R_2$ or not.

Setting $(X',Y',Z',W')=\left(q_1,\frac{p_1}{p_2},q_2,\frac{1}{p_2} \right)$, we will see the following relations;

Holomorphy condition $R_1$ $\Longleftrightarrow$ Accessible sing. $(X',Y',Z',W')=(X',0,0,0)$, \ Local index:$\begin{pmatrix}
\begin{matrix}
0 \\
1 \\
2 \\
1
\end{matrix}
\end{pmatrix}$, 

Holomorphy condition $R_4$ $\Longleftrightarrow$ Accessible sing. $(X',Y',Z',W')=(0,Y',0,0)$, \ Local index:$\begin{pmatrix}
\begin{matrix}
1 \\
0 \\
2 \\
1
\end{matrix}
\end{pmatrix}$.

It is also still an open question whether in Lax pair theory we can distinguish Holomorphy condition $R_1$ from Holomorphy condition $R_4$ or not.

}

\section{Appendix F: Holomorphy of Garnier system}

In this appendix, we remark some holomorphy conditions of the system \eqref{eq:Gar}.

{\bf Hamiltonians $H_{24}^{(1)}=r_{24}(H_1), \ H_{24}^{(2)}=r_{24}(H_2-p_1), \ r_{24}:x=-((q_1-s)p_1-\alpha_4)p_1, \ y=\frac{1}{p_1}, \ z=-(q_2p_2-\alpha_2)p_2, \ w=\frac{1}{p_2}$}
\begin{align*}
&r_0^{24}:x_0=-(q_1p_1+q_2 p_2-(\alpha_0+\alpha_2+\alpha_4))p_1, \ y_0=\frac{1}{p_1}, \ z_0=q_2 p_1, \ w_0=\frac{p_2}{p_1}, \\
&r_1^{24}:x_1=q_1, \ y_1=p_1, \ z_1=q_2+\frac{\alpha_1-\alpha_2}{p_2}+\frac{1}{p_2^2}, \ w_1=p_2, \\
&r_2^{24}:x_2=q_1, \ y_2=p_1, \ z_2=-(q_2 p_2-\alpha_2)p_2 \ w_2=\frac{1}{p_2}, \\
&r_3^{24}:x_3=q_1+q_2+\frac{2q_2p_2+\alpha_3-\alpha_2-\alpha_4}{p_1}-\frac{s}{p_1^2}, \ y_3=p_1, \ z_3=q_2p_1^2, \ w_3=\frac{p_2+p_1}{p_1^2}, \\
&r_4^{24}:x_4=-(q_1 p_1-\alpha_4)p_1, \ y_4=\frac{1}{p_1}, \ z_4=q_2, \ w_4=p_2, \\
&r_5^{24}:x_5=q_1+\frac{\alpha_5-\alpha_4}{p_1}+\frac{t-s}{p_1^2}, \ y_5=p_1, \ z_5=q_2, \ w_5=p_2,
\end{align*}
where $r_5^{24} \left(H_{24}^{(1)}-\frac{1}{p_1} \right), \ r_3^{24} \left(H_{24}^{(2)}+\frac{1}{p_1} \right), \ r_5^{24} \left(H_{24}^{(2)}+\frac{1}{p_1} \right)$.

\hspace{0.5cm}

{\bf Hamiltonians $H_{23}^{(1)}=r_{23}(H_1), \ H_{23}^{(2)}=r_{23}(H_2), \ r_{23}:x=-((q_1-q_2)p_1-\alpha_3)p_1, \ y=\frac{1}{p_1}, \ z=-(q_2(p_2+p_1)-\alpha_2)(p_2+p_1), \ w=\frac{1}{p_2+p_1}$}
\begin{align*}
&r_0^{23}:x_0=-(q_1p_1+q_2 p_2-(\alpha_0+\alpha_2+\alpha_3))p_1, \ y_0=\frac{1}{p_1}, \ z_0=q_2 p_1, \ w_0=\frac{p_2}{p_1}, \\
&r_1^{23}:x_1=q_1, \ y_1=p_1, \ z_1=q_2+\frac{\alpha_1-\alpha_2}{p_2}+\frac{1}{p_2^2}, \ w_1=p_2, \\
&r_2^{23}:x_2=q_1, \ y_2=p_1, \ z_2=-(q_2 p_2-\alpha_2)p_2 \ w_2=\frac{1}{p_2}, \\
&r_3^{23}:x_3=-(q_1 p_1-\alpha_3)p_1, \ y_3=\frac{1}{p_1}, \ z_3=q_2, \ w_3=p_2, \\
&r_4^{23}:x_4=q_1-q_2+\frac{2q_2p_2+\alpha_4-\alpha_2-\alpha_3}{p_1}+\frac{s}{p_1^2}, \ y_5=p_1, \ z_5=q_2p_1^2, \ w_5=\frac{p_2-p_1}{p_1^2}, \\
&r_5^{23}:x_5=q_1-q_2+\frac{2q_2p_2+\alpha_5-\alpha_2-\alpha_3}{p_1}+\frac{t}{p_1^2}, \ y_5=p_1, \ z_5=q_2p_1^2, \ w_5=\frac{p_2-p_1}{p_1^2},
\end{align*}
where $r_5^{23} \left(H_{23}^{(1)}-\frac{1}{p_1} \right), \ r_4^{23} \left(H_{24}^{(2)}-\frac{1}{p_1} \right)$.


\begin{thebibliography}{99}

\bibitem[1]{P1} P. Painlev\'e, {\em M\'emoire sur les \'equations diff\'erentielles dont l'int\'egrale g\'en\'erale est uniforme}, Bull. Soci\'et\'e Math\'ematique de France. {\bf 28} (1900),  201--261. 


\bibitem[2]{P2}  P. Painlev\'e, {\em Sur les \'equations diff\'erentielles du second ordre et d'ordre sup\'erieur dont l'int\'egrale est uniforme}, Acta Math. {\bf 25} (1902), 1--85. 


\bibitem[3]{BG} B. Gambier, {\em Sur les \'equations diff\'erentielles du second ordre et du premier degr\'e dont l'int\'egrale g\'en\'erale est \`a points critiques fixes}, Acta Math. {\bf 33} (1910), 1--55.

\bibitem[4]{Gar1} R. Garnier, 
{\em Sur des equations differentielles du troisieme ordre dont l'integrale generale est uniforme et sur une classe d'equations nouvelles d'ordre superieur dont l'integrale generale a ses points critiques fixes}, Ann. Sci. Ecole Norm. Sup. {\bf 29} (1912), 1--126.




\bibitem[5]{cos1} C. M. Cosgrove, 
{\em Higher-order Painlev\'e equations in the polynomial class I. Bureau symbol P2}, Stud. Appl. Math. {\bf 104}, (2000), 1--65.

\bibitem[6]{cos2} C. M. Cosgrove,
{\em Higher order Painlev\'e equations in the polynomial class II, Bureau symbol P1}, Studies in Applied Mathematics. {\bf 116} (2006).


\bibitem[7]{chazy} J. Chazy, 
{\em Sur les \'equations diff\'erentielles du trousi\'eme ordre et d'ordre sup\'erieur dont l'int\'egrale a ses points critiques fixes}, 
Acta Math. {\bf 34} (1911), 317--385.


\bibitem[8]{FB} F. Bureau, 
{\em Integration of some nonlinear systems of ordinary differential equations}, 
Annali di Matematica. {\bf 94} (1972), 345--359. 



\bibitem[9]{Lax} P. D. Lax, 
{\em Integrals of nonlinear equations of evolution and solitary waves}, Commun. Pure Appl. Math. {\bf 21} (1968), 467--490.


\bibitem[10]{N1} M. Noumi and Y. Yamada, 
        {\em Higher order Painlev\'e equations of type $A_l^{(1)}$}, Funkcial. Ekvac. {\bf 41} (1998), 483--503.

\bibitem[11]{Ta} N. Tahara, {\em An augmentation of the phase space of the system of type $A_4^{(1)}$}, Kyushu J. Math. {\bf 58} (2004), 393--425.

\bibitem[12]{Ince} E. L. Ince, 
{\em Ordinary differential equations}, Dover Publications, New York, (1956). 


\bibitem[13]{sasa2} Y. Sasano, {\em Higher order Painlev\'e equations of type ${D_l}^{(1)}$}, RIMS Kokyuroku {\bf 1473} (2006), 143--163. 

\bibitem[14]{sasa3} Y. Sasano, {\em Coupled Painlev\'e VI systems in dimension four with affine Weyl group symmetry of type $D_6^{(1)}$, II}, RIMS Kokyuroku Bessatsu. {\bf B5} (2008), 137--152.


\bibitem[15]{O3} K. Okamoto, {\em Sur les 
feuilletages associ\'es aux \'equations du second ordre \`a points critiques fixes de P. Painlev\'e, Espaces des conditions initiales}, Japan. 
J. Math. {\bf 5} (1979), 1--79. 

\bibitem[16]{T1} T. Shioda and K. Takano, {\em On some Hamiltonian structures of Painlev\'e systems I}, Funkcial. Ekvac. {\bf 40} (1997), 271--291.

\bibitem[17]{MMT} T. Matano, A. Matumiya and K. Takano, 
{\em On some Hamiltonian structures of Painlev\'e systems, II}, 
J. Math. Soc. Japan {\bf 51} (1999), 843--866.


\bibitem[18]{N2} M. Noumi and Y. Yamada,
        {\em Affine Weyl Groups, Discrete Dynamical Systems and Painlev\`e Equations}, Comm Math Phys {\bf 199} (1998), 281--295.

\bibitem[19]{Fuji} K. Fuji and T. Suzuki, {\em Higher order Painlev\'e system of type $D_{2n+2}^{(1)}$ arising from integrable hierarchy}, Int. Math. Res. Not. {\bf 1}(2008), Art.ID rnm129.

\bibitem[20]{Oshima} T. Oshima, {\em Classification of Fuchsian systems and their connection problem},  RIMS Kokyuroku Bessatsu {\bf B37} (2013), 163--192.


\bibitem[21]{Oshima2} K. Hiroe and T. Oshima, {\em A classification of roots of symmetric Kac-Moody root systems and its application}, Symmetries, Integral Systems and Representations, Springer Proceedings of Mathematics and Statics {\bf 40 }(2012), 195--241.

\bibitem[22]{KNS} Hiroshi Kawakami, Akane Nakamura, and Hidetaka Sakai, {\em Degeneration scheme of 4-dimensional Painlev\'e-type equations}, to appear in Contemporary Mathematics.

\bibitem[23]{Fuji2} K. Fuji, {\em private communication}.


\bibitem[24]{K} H. Kimura, {\em Uniform foliation associated with the Hamiltonian system ${\mathcal H}_{n}$}, Ann. Scuola Norm. Sup. Pisa Cl. Sci. (4) {\bf 20} (1993), no. 1, 1--60.


\bibitem[25]{KOka} H. Kimura and K. Okamoto, {\em On the polynomial Hamiltonian structure of the Garnier systems}, J. Math. Pures Appl. {\bf 63} (1984), 129--146.

\bibitem[26]{S} T. Suzuki, {\em Affine Weyl group symmetry of the Garnier system}, Funkcial. Ekvac. {\bf 48} (2005), 203--230.

\bibitem[27]{oka4} K. Okamoto, {\em Isomonodromic deformations and Painlev\'e equations, and the Garnier system}, J. Fac. Sci. Univ. Tokyo, Sect. IA Math., {\bf 33} (1986), 575--618.

\bibitem[27]{sasa7} Y. Sasano, {\em Coupled Painlev\'e IV systems in dimension four},
Kumamoto J. Math. {\bf 20} (2007) 13--31.

\bibitem[28]{Yamada} Y. Yamada and Y. Sasano, {\em Symmetry and holomorphy of Painlev\'e type systems}, RIMS Kokyuroku Bessatsu. {\bf B2} (2007) 215--225.

\bibitem[29]{JM1} M. Jimbo and T. Miwa, 
{\em Monodromy preserving deformation of linear ordinary differential equations with rational coefficients.II}, 
Physica, {\bf 2D} (1981), 407--448.

\bibitem[30]{Tsuda1} T. Tsuda, {\em Birational Symmetries, Hirota Bilinear Forms and Special Solutions of the Garnier Systems in 2-variables}, J. Math .Sci. Univ. Tokyo. {\bf 10} (2003), 355--371.

\bibitem[31]{Tsuda2} T. Tsuda, {\em Rational Solutions of the Garnier System in Terms of Schur Polynomials}, IMRN. {\bf 43} (2003), 2341--2358.

\bibitem[32]{Tsuda3} T. Tsuda, {\em Universal characters and integrable systems}, PhD thesis. The University of Tokyo, (2003).

\bibitem[33]{Tsuda4} T. Tsuda, {\em Toda equation and special polynomials associated with the Garnier system}, Advances in Mathematics. {\bf 206} (2006), 

\bibitem[34]{Candelas} Philip CANDELAS, Xenia C. DE LA OSSA, Paul S. GREEN and Linda PARKES, {\em A PAIR OF CALABI-YAU MANIFOLDS AS AN EXACTLY SOLUBLE
SUPERCONFORMAL THEORY}, Nuclear Physics. {\bf B359} (1991) 21--74 North-Holland.

\bibitem[35]{FS1} Kenta Fuji ,Takao Suzuki, {\em Drinfeld-Sokolov hierarchies of type A and fourth order Painlev\'e systems}, Funkcial. Ekvac {\bf 53} (2010), 143--167.

\bibitem[36]{sasa8} Y. Sasano, {\em Painlev\'e scheme}, arXiv:0709.0597.

\end{thebibliography}
\end{document}